\documentclass[11pt,leqno]{article}
\normalfont
\renewcommand{\baselinestretch}{1.24}

\usepackage{graphicx}
\usepackage{amsfonts}
\usepackage{mathrsfs}
\usepackage{rotating}
\usepackage{mathdots}
\usepackage{hieroglf}
\usepackage{MnSymbol}
\usepackage[dvipsnames]{xcolor}
\usepackage{extarrows}

\newcommand{\Intro}{1}
\newcommand{\Background}{2}
\newcommand{\MoveMinGame}{3}
\newcommand{\MixedMiddleSwitch}{4}
\newcommand{\Closing}{5}

\newcommand{\TopoBalanceFigure}{Figure 2.1}
\newcommand{\MountainValleyLemma}{Lemma 2.1}
\newcommand{\ModularLatticeProp}{Proposition 2.2}
\newcommand{\JcolorFigure}{Figure 2.2}
\newcommand{\TechnicalLemma}{Proposition 2.3}
\newcommand{\ColorComponentTheorem}{Theorem 2.4}
\newcommand{\DilworthCorollary}{Corollary 2.5}

\newcommand{\MoveMinGameTheorem}{Theorem 3.1}
\newcommand{\MoveMinSectionRemarks}{Remarks 3.2}

\newcommand{\TypeBLatticeIsoProp}{Proposition 4.1}
\newcommand{\MixedMiddleSwitchProp}{Proposition 4.2} 
\newcommand{\BFiveFigure}{Figure 4.1}

\newcommand{\TypeBIsoAndMixedmiddleProps}{Propositions 4.1 and 4.2}
\newcommand{\ChallengeFigure}{Figure 5.1}
\newcommand{\CatalanLatticeFigure}{Figure 5.2}
\newcommand{\CatalanIsoProposition}{Proposition 5.1}

\newfont{\myscbolditalics}{ecoc0500 at 11pt}

\newfont{\mybolditalics}{ecbi0500 at 11pt}

\newfont{\mysmallbolditalics}{ecbi0500 at 9pt}

\newcommand{\myspecialchar}[1]{\mbox{\myscbolditalics #1}}

\newfont{\mysmallscbolditalics}{ecoc0500 at 8pt}

\newfont{\eulercursive}{eurm10 at 11pt}
\newcommand{\myd}{\mbox{\eulercursive d}}
\newcommand{\mya}{\mbox{\eulercursive a}}
\newcommand{\mysmallerd}{\mbox{\smalleulercursive d}}
\newcommand{\mysmallera}{\mbox{\smalleulercursive a}}
\newcommand{\mys}{\mbox{\eulercursive s}}
\newcommand{\myt}{\mbox{\eulercursive t}}
\newcommand{\myx}{\mbox{\eulercursive x}}
\newcommand{\myy}{\mbox{\eulercursive y}}

\newcommand{\myl}{\mbox{\eulercursive `}}
\newfont{\smalleulercursive}{eurm10 at 9pt}

\newfont{\smallereulercursive}{eurm10 at 7pt}

\newfont{\largereulercursive}{eurm10 at 14pt}
\newcommand{\mylargers}{\mbox{\largereulercursive s}}
\newcommand{\mylargert}{\mbox{\largereulercursive t}}

\newfont{\myslantcyrillic}{wncyi10 at 11pt}

\newcommand{\QED}{\raisebox{0.5mm}{\fbox{\rule{0mm}{1.5mm}\ }}}

\newcounter{myfn}[page]
\renewcommand{\thefootnote}{\fnsymbol{footnote}}

\newcounter{rone}
\newcounter{rtwo}
\newcounter{rthree}
\newcounter{rfour}
\newcounter{rfive}
\newcounter{rsix}
\newcounter{rseven}
\newcommand{\myB}{\mbox{\sffamily B}}

\newcommand{\aelt}{\mathbf{a}}

\newcommand{\melt}{\mathbf{m}} \newcommand{\nelt}{\mathbf{n}}
 
 \newcommand{\relt}{\mathbf{r}}
\newcommand{\selt}{\mathbf{s}} \newcommand{\telt}{\mathbf{t}}
\newcommand{\uelt}{\mathbf{u}} \newcommand{\velt}{\mathbf{v}}
\newcommand{\welt}{\mathbf{w}} \newcommand{\xelt}{\mathbf{x}}
\newcommand{\yelt}{\mathbf{y}} \newcommand{\zelt}{\mathbf{z}}

\newcommand{\ecolor}{\mbox{\sffamily edge{\_}color}}
\newcommand{\vcolor}{\mbox{\sffamily vertex{\_}color}}
\newcommand{\dist}{\mbox{\sffamily dist}}
\newcommand{\comp}{\mbox{\sffamily comp}}

\newcommand{\myarrow}[1]{\stackrel{#1}{\rightarrow}}

\newcommand{\mylongarrow}[1]{\stackrel{#1}{\longrightarrow}}
\newcommand{\mylongbackarrow}[1]{\stackrel{\ #1}{\longleftarrow}}

\newcommand{\eqmulti}{\xlongequal{\text{multiset}}}%

\newcommand{\pathlength}{\mbox{\sffamily path{\_}length}}

\newcommand{\change}{\mbox{\sffamily change}}

\newcommand{\CatalanVertexRankZero}[2]{
\setlength{\unitlength}{0.2cm}
\begin{picture}(0,0)
\put(0,0){\circle*{0.5}} 
\put(0,-1.5){\put(#1,#2){\begin{picture}(0,0)
\put(0,0){\BallotZero{0.175cm}}
\end{picture}}}
\end{picture}
}
\newcommand{\CatalanVertexRankOne}[2]{
\setlength{\unitlength}{0.2cm}
\begin{picture}(0,0)
\put(0,0){\circle*{0.5}} 
\put(0,-1.5){\put(#1,#2){\begin{picture}(0,0)
\put(0,0){\BallotOne{0.175cm}}
\end{picture}}}
\end{picture}
}
\newcommand{\CatalanVertexRankTwoNumOne}[2]{
\setlength{\unitlength}{0.2cm}
\begin{picture}(0,0)
\put(0,0){\circle*{0.5}} 
\put(0,-1.5){\put(#1,#2){\begin{picture}(0,0)
\put(0,0){\BallotTwoOne{0.175cm}}
\end{picture}}}
\end{picture}
}
\newcommand{\CatalanVertexRankTwoNumTwo}[2]{
\setlength{\unitlength}{0.2cm}
\begin{picture}(0,0)
\put(0,0){\circle*{0.5}} 
\put(0,-1.5){\put(#1,#2){\begin{picture}(0,0)
\put(0,0){\BallotTwoTwo{0.175cm}}
\end{picture}}}
\end{picture}
}
\newcommand{\CatalanVertexRankThreeNumOne}[2]{
\setlength{\unitlength}{0.2cm}
\begin{picture}(0,0)
\put(0,0){\circle*{0.5}} 
\put(0,-1.5){\put(#1,#2){\begin{picture}(0,0)
\put(0,0){\BallotThreeOne{0.175cm}}
\end{picture}}}
\end{picture}
}
\newcommand{\CatalanVertexRankThreeNumTwo}[2]{
\setlength{\unitlength}{0.2cm}
\begin{picture}(0,0)
\put(0,0){\circle*{0.5}} 
\put(0,-1.5){\put(#1,#2){\begin{picture}(0,0)
\put(0,0){\BallotThreeTwo{0.175cm}}
\end{picture}}}
\end{picture}
}
\newcommand{\CatalanVertexRankThreeNumThree}[2]{
\setlength{\unitlength}{0.2cm}
\begin{picture}(0,0)
\put(0,0){\circle*{0.5}} 
\put(0,-1.5){\put(#1,#2){\begin{picture}(0,0)
\put(0,0){\BallotThreeThree{0.175cm}}
\end{picture}}}
\end{picture}
}
\newcommand{\CatalanVertexRankFourNumOne}[2]{
\setlength{\unitlength}{0.2cm}
\begin{picture}(0,0)
\put(0,0){\circle*{0.5}} 
\put(0,-1.5){\put(#1,#2){\begin{picture}(0,0)
\put(0,0){\BallotFourOne{0.175cm}}
\end{picture}}}
\end{picture}
}
\newcommand{\CatalanVertexRankFourNumTwo}[2]{
\setlength{\unitlength}{0.2cm}
\begin{picture}(0,0)
\put(0,0){\circle*{0.5}} 
\put(0,-1.5){\put(#1,#2){\begin{picture}(0,0)
\put(0,0){\BallotFourTwo{0.175cm}}
\end{picture}}}
\end{picture}
}
\newcommand{\CatalanVertexRankFourNumThree}[2]{
\setlength{\unitlength}{0.2cm}
\begin{picture}(0,0)
\put(0,0){\circle*{0.5}} 
\put(0,-1.5){\put(#1,#2){\begin{picture}(0,0)
\put(0,0){\BallotFourThree{0.175cm}}
\end{picture}}}
\end{picture}
}
\newcommand{\CatalanVertexRankFiveNumOne}[2]{
\setlength{\unitlength}{0.2cm}
\begin{picture}(0,0)
\put(0,0){\circle*{0.5}} 
\put(0,-1.5){\put(#1,#2){\begin{picture}(0,0)
\put(0,0){\BallotFiveOne{0.175cm}}
\end{picture}}}
\end{picture}
}
\newcommand{\CatalanVertexRankFiveNumTwo}[2]{
\setlength{\unitlength}{0.2cm}
\begin{picture}(0,0)
\put(0,0){\circle*{0.5}} 
\put(0,-1.5){\put(#1,#2){\begin{picture}(0,0)
\put(0,0){\BallotFiveTwo{0.175cm}}
\end{picture}}}
\end{picture}
}
\newcommand{\CatalanVertexRankFiveNumThree}[2]{
\setlength{\unitlength}{0.2cm}
\begin{picture}(0,0)
\put(0,0){\circle*{0.5}} 
\put(0,-1.5){\put(#1,#2){\begin{picture}(0,0)
\put(0,0){\BallotFiveThree{0.175cm}}
\end{picture}}}
\end{picture}
}
\newcommand{\CatalanVertexRankSix}[2]{
\setlength{\unitlength}{0.2cm}
\begin{picture}(0,0)
\put(0,0){\circle*{0.5}} 
\put(0,-1.5){\put(#1,#2){\begin{picture}(0,0)
\put(0,0){\BallotSix{0.175cm}}
\end{picture}}}
\end{picture}
}

\newcommand{\VertexForBnLattice}[7]{
\setlength{\unitlength}{1.325cm}
\begin{picture}(0,0)
\put(-0,0){\circle*{0.12}} 
\put(0.25,0.05){\put(#6,#7){\tiny #1#2#3#4#5}}
\end{picture}
}

\newcommand{\VertexForBnLatticeToo}[7]{
\setlength{\unitlength}{1.325cm}
\begin{picture}(0,0)
\put(0.25,-0.1){\put(#6,#7){\textcolor{Gray}{\tiny #1#2#3#4#5}}}
\end{picture}
}

\newcommand{\NEEdgeForBnLattice}[4]{
\setlength{\unitlength}{1.325cm}
\begin{picture}(0,0)
\thicklines
\textcolor{#1}{
\put(0,0){\line(3,2){1}}
\put(0.4,0.3){\put(#3,#4){\small \em #2}} 
}
\end{picture}
}

\newcommand{\NWEdgeForBnLattice}[4]{
\setlength{\unitlength}{1.325cm}
\begin{picture}(0,0)
\thicklines
\textcolor{#1}{
\put(0,0){\line(-3,2){1}}
\put(-0.525,0.3){\put(#3,#4){\small \em #2}} 
}
\end{picture}
}

\newcommand{\VerticalEdgeForBnLattice}[4]{
\setlength{\unitlength}{1.325cm}
\begin{picture}(0,0)
\thicklines
\textcolor{#1}{
\put(0,0){\line(0,1){0.6667}}
\put(-0.05,0.3){\put(#3,#4){\small \em #2}} 
}
\end{picture}
}

\newcommand{\BallotZero}[1]{
\setlength{\unitlength}{#1}
\begin{picture}(3,3)
\put(0,0){\line(0,1){3}} \put(1,0){\line(0,1){3}} \put(2,1){\line(0,1){2}} \put(3,2){\line(0,1){1}}
\put(0,0){\line(1,0){1}} \put(0,1){\line(1,0){2}} \put(0,2){\line(1,0){3}} \put(0,3){\line(1,0){3}}
\end{picture}}

\newcommand{\BallotOne}[1]{
\setlength{\unitlength}{#1}
\begin{picture}(3,3)
\thicklines
\multiput(0.05,2)(0.05,0){19}{\textcolor{gray}{\line(0,1){1}}}
\thinlines
\put(0,0){\line(0,1){3}} \put(1,0){\line(0,1){3}} \put(2,1){\line(0,1){2}} \put(3,2){\line(0,1){1}}
\put(0,0){\line(1,0){1}} \put(0,1){\line(1,0){2}} \put(0,2){\line(1,0){3}} \put(0,3){\line(1,0){3}}
\end{picture}}

\newcommand{\BallotTwoOne}[1]{
\setlength{\unitlength}{#1}
\begin{picture}(3,3)
\thicklines
\multiput(0.05,1)(0.05,0){19}{\textcolor{gray}{\line(0,1){2}}}
\thinlines
\put(0,0){\line(0,1){3}} \put(1,0){\line(0,1){3}} \put(2,1){\line(0,1){2}} \put(3,2){\line(0,1){1}}
\put(0,0){\line(1,0){1}} \put(0,1){\line(1,0){2}} \put(0,2){\line(1,0){3}} \put(0,3){\line(1,0){3}}
\end{picture}}

\newcommand{\BallotTwoTwo}[1]{
\setlength{\unitlength}{#1}
\begin{picture}(3,3)
\thicklines
\multiput(0.05,2)(0.05,0){39}{\textcolor{gray}{\line(0,1){1}}}
\thinlines
\put(0,0){\line(0,1){3}} \put(1,0){\line(0,1){3}} \put(2,1){\line(0,1){2}} \put(3,2){\line(0,1){1}}
\put(0,0){\line(1,0){1}} \put(0,1){\line(1,0){2}} \put(0,2){\line(1,0){3}} \put(0,3){\line(1,0){3}}
\end{picture}}

\newcommand{\BallotThreeOne}[1]{
\setlength{\unitlength}{#1}
\begin{picture}(3,3)
\thicklines
\multiput(0.05,0)(0.05,0){19}{\textcolor{gray}{\line(0,1){3}}}
\thinlines
\put(0,0){\line(0,1){3}} \put(1,0){\line(0,1){3}} \put(2,1){\line(0,1){2}} \put(3,2){\line(0,1){1}}
\put(0,0){\line(1,0){1}} \put(0,1){\line(1,0){2}} \put(0,2){\line(1,0){3}} \put(0,3){\line(1,0){3}}
\end{picture}}

\newcommand{\BallotThreeTwo}[1]{
\setlength{\unitlength}{#1}
\begin{picture}(3,3)
\thicklines
\multiput(0.05,1)(0.05,0){19}{\textcolor{gray}{\line(0,1){2}}}
\multiput(1.05,2)(0.05,0){19}{\textcolor{gray}{\line(0,1){1}}}
\thinlines
\put(0,0){\line(0,1){3}} \put(1,0){\line(0,1){3}} \put(2,1){\line(0,1){2}} \put(3,2){\line(0,1){1}}
\put(0,0){\line(1,0){1}} \put(0,1){\line(1,0){2}} \put(0,2){\line(1,0){3}} \put(0,3){\line(1,0){3}}
\end{picture}}

\newcommand{\BallotThreeThree}[1]{
\setlength{\unitlength}{#1}
\begin{picture}(3,3)
\thicklines
\multiput(0.05,2)(0.05,0){59}{\textcolor{gray}{\line(0,1){1}}}
\thinlines
\put(0,0){\line(0,1){3}} \put(1,0){\line(0,1){3}} \put(2,1){\line(0,1){2}} \put(3,2){\line(0,1){1}}
\put(0,0){\line(1,0){1}} \put(0,1){\line(1,0){2}} \put(0,2){\line(1,0){3}} \put(0,3){\line(1,0){3}}
\end{picture}}

\newcommand{\BallotFourOne}[1]{
\setlength{\unitlength}{#1}
\begin{picture}(3,3)
\thicklines
\multiput(0.05,0)(0.05,0){19}{\textcolor{gray}{\line(0,1){3}}}
\multiput(1.05,2)(0.05,0){19}{\textcolor{gray}{\line(0,1){1}}}
\thinlines
\put(0,0){\line(0,1){3}} \put(1,0){\line(0,1){3}} \put(2,1){\line(0,1){2}} \put(3,2){\line(0,1){1}}
\put(0,0){\line(1,0){1}} \put(0,1){\line(1,0){2}} \put(0,2){\line(1,0){3}} \put(0,3){\line(1,0){3}}
\end{picture}}

\newcommand{\BallotFourTwo}[1]{
\setlength{\unitlength}{#1}
\begin{picture}(3,3)
\thicklines
\multiput(0.05,1)(0.05,0){39}{\textcolor{gray}{\line(0,1){2}}}
\thinlines
\put(0,0){\line(0,1){3}} \put(1,0){\line(0,1){3}} \put(2,1){\line(0,1){2}} \put(3,2){\line(0,1){1}}
\put(0,0){\line(1,0){1}} \put(0,1){\line(1,0){2}} \put(0,2){\line(1,0){3}} \put(0,3){\line(1,0){3}}
\end{picture}}

\newcommand{\BallotFourThree}[1]{
\setlength{\unitlength}{#1}
\begin{picture}(3,3)
\thicklines
\multiput(0.05,2)(0.05,0){59}{\textcolor{gray}{\line(0,1){1}}}
\multiput(0.05,1)(0.05,0){19}{\textcolor{gray}{\line(0,1){1}}}
\thinlines
\put(0,0){\line(0,1){3}} \put(1,0){\line(0,1){3}} \put(2,1){\line(0,1){2}} \put(3,2){\line(0,1){1}}
\put(0,0){\line(1,0){1}} \put(0,1){\line(1,0){2}} \put(0,2){\line(1,0){3}} \put(0,3){\line(1,0){3}}
\end{picture}}

\newcommand{\BallotFiveOne}[1]{
\setlength{\unitlength}{#1}
\begin{picture}(3,3)
\thicklines
\multiput(0.05,0)(0.05,0){19}{\textcolor{gray}{\line(0,1){3}}}
\multiput(1.05,1)(0.05,0){19}{\textcolor{gray}{\line(0,1){2}}}
\thinlines
\put(0,0){\line(0,1){3}} \put(1,0){\line(0,1){3}} \put(2,1){\line(0,1){2}} \put(3,2){\line(0,1){1}}
\put(0,0){\line(1,0){1}} \put(0,1){\line(1,0){2}} \put(0,2){\line(1,0){3}} \put(0,3){\line(1,0){3}}
\end{picture}}

\newcommand{\BallotFiveTwo}[1]{
\setlength{\unitlength}{#1}
\begin{picture}(3,3)
\thicklines
\multiput(0.05,0)(0.05,0){19}{\textcolor{gray}{\line(0,1){3}}}
\multiput(1.05,2)(0.05,0){39}{\textcolor{gray}{\line(0,1){1}}}
\thinlines
\put(0,0){\line(0,1){3}} \put(1,0){\line(0,1){3}} \put(2,1){\line(0,1){2}} \put(3,2){\line(0,1){1}}
\put(0,0){\line(1,0){1}} \put(0,1){\line(1,0){2}} \put(0,2){\line(1,0){3}} \put(0,3){\line(1,0){3}}
\end{picture}}

\newcommand{\BallotFiveThree}[1]{
\setlength{\unitlength}{#1}
\begin{picture}(3,3)
\thicklines
\multiput(0.05,2)(0.05,0){59}{\textcolor{gray}{\line(0,1){1}}}
\multiput(0.05,1)(0.05,0){39}{\textcolor{gray}{\line(0,1){1}}}
\thinlines
\put(0,0){\line(0,1){3}} \put(1,0){\line(0,1){3}} \put(2,1){\line(0,1){2}} \put(3,2){\line(0,1){1}}
\put(0,0){\line(1,0){1}} \put(0,1){\line(1,0){2}} \put(0,2){\line(1,0){3}} \put(0,3){\line(1,0){3}}
\end{picture}}

\newcommand{\BallotSix}[1]{
\setlength{\unitlength}{#1}
\begin{picture}(3,3)
\thicklines
\multiput(0.05,2)(0.05,0){59}{\textcolor{gray}{\line(0,1){1}}}
\multiput(0.05,1)(0.05,0){39}{\textcolor{gray}{\line(0,1){1}}}
\multiput(0.05,0)(0.05,0){19}{\textcolor{gray}{\line(0,1){1}}}
\thinlines
\put(0,0){\line(0,1){3}} \put(1,0){\line(0,1){3}} \put(2,1){\line(0,1){2}} \put(3,2){\line(0,1){1}}
\put(0,0){\line(1,0){1}} \put(0,1){\line(1,0){2}} \put(0,2){\line(1,0){3}} \put(0,3){\line(1,0){3}}
\end{picture}}

\begin{document}
\pagenumbering{arabic}
\thispagestyle{empty}%
\vspace*{-0.7in}
\noindent
\hfill {\scriptsize December 3, 2023}

\begin{center}
{\large \bf Move-minimizing puzzles and diamond-colored modular $\!/\!$ distributive lattices} 

\vspace*{0.05in}
\renewcommand{\thefootnote}{1}
Robert G.\ Donnelly\footnote{Department of Mathematics and Statistics, Murray State
University, Murray, KY 42071\\ 
\hspace*{0.25in}Email: {\tt rob.donnelly@murraystate.edu}}$\!$, 
\renewcommand{\thefootnote}{2}
\hspace*{-0.05in}Elizabeth A.\ Donovan\footnote{Department of Mathematics and Statistics, Murray State University, Murray, KY 42071\\ 
\hspace*{0.25in}Email: {\tt edonovan@murraystate.edu},\ Fax: 1-270-809-2314}$\!$, 
\renewcommand{\thefootnote}{3} 
\hspace*{-0.05in}Molly W.\ Dunkum\footnote{Department of Mathematics, Western Kentucky University, Bowling Green, KY 42101\\ 
\hspace*{0.25in}Corresponding author,\ Email: {\tt molly.dunkum@wku.edu}}$\!$, and 
\renewcommand{\thefootnote}{4}
\hspace*{-0.05in}Timothy A.\ Schroeder\footnote{Department of Mathematics and Statistics, Murray State University, Murray, KY 42071\\ 
\hspace*{0.25in}Email: {\tt tschroeder@murraystate.edu}}

\end{center} 

\begin{abstract}
The move-minimizing puzzles presented here are certain types of one-player combinatorial games that are shown to have explicit solutions whenever they can be encoded in a certain way as diamond-colored modular or distributive lattices. 
Our work here is founded in a new interpretation of some routine and elementary order-theoretic combinatorics.

\begin{center}

\vspace*{0.1in}

\end{center} 
\end{abstract}

\vspace{1ex} 
\noindent 
{\Large \bf \Intro\hspace*{0.15in} Introduction}

\noindent
Imagine that a game board consists of a set of switches that can be in one of two states, say `ON' or `OFF'.  
A {\em switching puzzle} on such a board might ask a player to maneuver from one ON/OFF arrangement of switches to another arrangement by changing the binary states of switches (i.e.\ by `toggling' switches from OFF to ON or vice-versa) in some sequence chosen by the player; most likely there will be a set of rules governing when a given switch can be toggled, the effects such toggling will have on near-by switches, etc.  
Such challenges are popular in electronic games; perhaps the most well-known such game is {\sl Lights Out} \cite{WikiLights}.  
Some further aspects of switching puzzles can be to decide whether it is even possible to maneuver from a given arrangement to some other arrangement and, if so, what minimum number of moves is required. 

For us, a {\em move-minimizing puzzle} is a combinatorial puzzle which consists of a set of {\em objects} (say, the possible arrangements of switches on a board) and a set of {\em legal moves} (which allow the player to transition from one object to another by, say, toggling switches) and which challenges the player to maneuver from a given object to another given object (perhaps in a minimum number of legal moves, if possible). 
Switching puzzles, as described in the previous paragraph, can be move-minimizing puzzles. 

It is natural to ask several questions concerning such a puzzle, such as: 
(1a) Is it always possible to maneuver from any  given object $\selt$ to any other given object $\telt$ using legal moves?  
If so, what is an explicit formula for the minimum number of legal moves required to move from $\selt$ to $\telt$? 
And, what is the maximum of all these `move-minimizing numbers' when all possible ordered pairs $(\selt,\telt)$ of objects are considered?  
(1b) When it is possible to maneuver from $\selt$ to $\telt$, can a sequence of legal moves from $\selt$ to $\telt$ be explicitly prescribed in such a way that it achieves this minimum number of uses?  
For our purposes, to answer these questions is to solve the given move-minimizing puzzle. 
One of our main results (\MoveMinGameTheorem) offers explicit answers to questions (1a) and (1b) for certain types of move-minimizing puzzles. 
Puzzles of this certain type are those that can be modelled by so-called `diamond-colored' modular or distributive lattices and therefore have the favorable properties developed in our other main result, \ColorComponentTheorem. 

The next set of questions might seem, at first, to be a trivial modification of questions (1a) and (1b) above.  
However, the structural advantages to the types of move-minimizing puzzles we study here will lead to perhaps-simpler-than-expected answers.  
(2a) If we restrict puzzle play by allowing only some subset $\mathcal{S}$ of legal moves to be used, under what conditions is it possible to maneuver from our given $\selt$ to our given $\telt$? 
If this is possible, what is an explicit formula for the minimum number of moves required under this restriction? 
Additionally, if we regard moves from the subset $\mathcal{S}$ as special but not required, then what is an explicit formula for the minimum number of times moves from $\mathcal{S}$ must be used in maneuvering from $\selt$ to $\telt$?  
(2b) If only moves from the given subset $\mathcal{S}$ of legal moves are allowed in maneuvering from $\selt$ to $\telt$, can a sequence of such moves be explicitly prescribed in such a way that it achieves the minimum number of moves? 
And, if moves from $\mathcal{S}$ are regarded as special but are not required in maneuvering from $\selt$ to $\telt$, can a sequence of such moves be explicitly prescribed in such a way that the minimum number of uses of the special moves is achieved? 
\MoveMinGameTheorem\ addresses questions (2a) and (2b) as well.  

In computer-science-oriented graph theory, such move-minimizing puzzles would be called `reconfiguration problems,' with our \underline{objects} playing the role of \underline{configurations}.  
Generally speaking, reconfiguration problems consider the scenario of moving from one given (feasible) configuration to another while maintaining feasibility of all intermediate configurations. 
For example, the challenge of repairing a network without interrupting service is essentially a reconfiguration problem. 
For a recent survey, see \cite{Nishimura}. 
Some reconfiguration problems are called `combinatorial games' and require two or more players (such as {\sl Checkers}); some are called `combinatorial puzzles' (also known as one-player combinatorial games, such as the well-known {\sl 15-puzzle}, which is a sliding puzzle); and some are `deterministic simulations' (or zero-player games, such as J.\ H.\ Conway's {\sl Game of Life}). 
Much of the literature on such reconfiguration games is devoted to analysis of their (time and space) complexity and related computational aspects, see for example \cite{Fraenkel}, \cite{HD}, or \cite{vdH}. 
All of the puzzles considered here are essentially {\sc All-Pairs Shortest Path} problems as described for example in \cite{CLRS}. 
Therefore, in principle, each puzzle can be solved in polynomial time by, say, the Floyd-(Ingerman-Roy-)Warshall algorithm. 
However, our purpose is to consider move-minimizing puzzles associated with a special type of graphs -- namely, covering digraphs of finite modular and distributive lattices -- in view of their distinctive order-theoretic features which afford explicit answers to our move-minimizing questions; for us, complexity is not a primary motivator. 

We now propose a switching puzzle. 
Let us fix an integer $n > 1$, and let $\mathcal{B}(n)$ be the set of all binary sequences of length $n$.  Given a binary sequence $(s_{1}, s_{2}, \ldots, s_{n-1}, s_{n})$, a {\em mixedmiddleswitch move} is the replacement of some sequence entry $s_{i}$ with $1-s_{i}$ subject to the following rules:  \fbox{$1 < i < n$} $s_{i}$ can be switched to $1-s_{i}$ only if $s_{i-1} \not= s_{i+1}$ (hence our idiosyncratic compound word `mixedmiddleswitch'); \fbox{$i=1$} $s_{1}$ can be switched to $1-s_{1}$ only if $s_{2} = 1$ (this can be viewed as a mixedmiddleswitch move if we regard our sequence as having a fixed $0^{\mbox{\tiny th}}$ entry $s_{0} := 0$); and \fbox{$i=n$} there are no restrictions on switching $s_{n}$ to $1-s_{n}$. 
The {\em Mixedmiddleswitch Puzzle} on $n$ switches is the move-minimizing puzzle whose objects are the binary sequences of $\mathcal{B}(n)$ and whose legal moves are the mixedmiddleswitch moves. 
At this point, we invite the reader to explore this puzzle for small values of $n$.  In particular we ask: In the Mixedmiddleswitch Puzzle on $n=5$ switches, can one maneuver from the binary sequence $(0,0,0,0,0)$ to the binary sequence $(0,1,0,1,0)$ using mixedmiddleswitch moves, and, if so, what minimum number of mixedmiddleswitch moves is required? 
This question is easily answered by applying \MixedMiddleSwitchProp\ below.  

We note, in passing, that an ulterior motive for our consideration of such puzzles is to find new models for certain algebraic-combinatorial objects, particularly representations of finite-dimensional semisimple Lie algebras over $\mathbb{C}$ and their companion Weyl group symmetric functions. 
Indeed, our Mixedmiddleswitch Puzzle is an example of such a connection, as our lattice-based solution to the puzzle concerns a classical family of diamond-colored distributive lattices we call the `type $\myB$ minuscule lattices' because of their affiliation with the minuscule representations of the odd orthogonal Lie algebras and the associated type $\myB$ Weyl group symmetric functions. 
Such algebraic connections are not formally considered here and are briefly taken up in the sequel paper on Domino Puzzles cf.\ \cite{BDDDS}; they will be developed most explicitly in \cite{BDD}. 

The main result of this paper is \MoveMinGameTheorem, which says generally how to solve those move-minimizing puzzles whose underlying graphs are covering digraphs for so-called diamond-colored modular/distributive lattices. 
We believe that our lattice-oriented approach to solving move-minimizing puzzles as framed in \MoveMinGameTheorem\ is new and might be of interest to a general audience readership. 
Our proof of \MoveMinGameTheorem\ relies on some standard order theory as well as some new\footnote{`New' in the sense of only have appeared in the first-named author's arXiv monograph \cite{DonDistributive}.} results concerning diamond-colored modular and distributive lattices that are developed in Section ~\Background\ (especially \ColorComponentTheorem). 
The Mixedmiddleswitch Puzzle, which is solved in Section~\MixedMiddleSwitch, will serve as our main example. 
In Section~\Closing, we present another example as a recreational challenge and as a further illustration of our lattice-theoretic perspective. 
The reader is strongly encouraged to play the move-minimizing games presented here by using small examples and to understand the various lattice isomorphisms included herein by working out these correspondences by hand for small cases.

\vspace{1ex} 
\noindent 
{\Large \bf \Background\hspace*{0.15in}Some background on digraphs, posets, and lattices}

\noindent 
To set up our main result (\MoveMinGameTheorem) and clarify our language/notation, we develop some preliminary concepts including some new results (\TechnicalLemma, \ColorComponentTheorem). 
The central combinatorial structures of this paper can ultimately be viewed as graphs.  
Generally speaking, any graph we work with will be finite and directed, with no loops and at most one edge between any two vertices, i.e.\ simple directed graphs.  
We will mainly work with partially ordered sets, thought of as directed graphs when identified with their covering digraphs. 
Most often, such a graph will either have its edges or its vertices `colored' by elements from some index set (usually a set of positive integers).  
Such coloring provides crucial information when we view these structures within the algebraic contexts that primarily motivate our interest. 
The conventions and notation we use concerning such structures largely follow the texts \cite{Aigner} and \cite{StanleyText} as well as the papers \cite{DonSupp}, \cite{ADLP}, and \cite{ADLMPPW}. 

Much of the work in this section is derived from the unpublished research monograph \cite{DonDistributive}. 
Part I of that document is a collection of elementary but foundational results relating to `diamond-colored' modular and distributive lattices intended to provide a browsable tutorial on the rudiments of this topic and to serve as a convenient single-source reference for readers interested in the supporting details. 
Most of those results are intuitive, are natural extensions of standard finite modular and distributive lattice concepts, and have antecedents in the literature in work of Proctor, Stembridge, and others (see e.g.\ \cite{PrEur}, \cite{PrGZ}, \cite{StemQuasi}, \cite{StemFC}, \cite{DonSupp}, \cite{ADLP}, \cite{ADLMPPW}). 
Of the results from \cite{DonDistributive} that are recapitulated here, \ColorComponentTheorem\ is the most efficacious.  
Other ideas and results -- such as the well-known claims made in \ModularLatticeProp, the {\bf mountain{\_}path{\_}construction} and {\bf valley{\_}path{\_}construction} algorithms which are the basis for most of the proofs of this section, and \MountainValleyLemma\ -- are straightforward and surely not novel. 

Suppose $R$ is a simple directed graph with vertex set $\mbox{\sffamily Vertices}(R)$ and directed edge set $\mbox{\sffamily Edges}(R)$.  
If $R$ is accompanied by a set $I$ of `colors' and a function $\ecolor_{R}: \mbox{\sffamily Edges}(R) \longrightarrow I$ (respectively $\vcolor_{R}: \mbox{\sffamily Vertices}(R) \longrightarrow I$), then we say $R$ is {\em edge-colored}  (resp.\ {\em vertex-colored}) {\em by the set} $I$. 
Denote by $R^{*}$ the {\em dual} edge-colored (resp., vertex-colored) digraph with $\mbox{\sffamily Vertices}(R^{*}) := \{\xelt^{*}\, |\, \xelt \in \mbox{\sffamily Vertices}(R)\}$, $\mbox{\sffamily Edges}(R^{*}) := \{\telt^{*} \rightarrow \selt^{*}\, |\, \selt \rightarrow \telt \in \mbox{\sffamily Edges}(R)\}$, and $\ecolor_{R^{*}}(\telt^{*} \rightarrow \selt^{*}) := \ecolor_{R}(\selt \rightarrow \telt)$ (resp., $\vcolor_{R^{*}}(\xelt^{*}) := \vcolor_{R}(\xelt)$). 
A {\em path} from $\selt$ to $\telt$ in $R$ is a sequence $\mathcal{P} = (\selt = \xelt_{0}, \xelt_{1}, \ldots, \xelt_{k} = \telt)$ such that 
for $j \in [1,k]_{\mathbb{Z}}$\footnote{For brevity, we denote by ${\mathcal{I}}_{\mathbb{Z}}$ the set of integers within a given real number interval $\mathcal{I}$.} we have either $\xelt_{j-1} \myarrow{i_j} \xelt_{j}$ or $\xelt_{j} \myarrow{i_j} \xelt_{j-1}$, where $(i_{j})_{j=1}^{k}$ is a sequence of colors from $I$.  
This path has {\em length} $k$, written $\pathlength(\mathcal{P})$, and we allow paths to have length $0$. 
Call $\mathcal{P}$ {\em simple} if each vertex of the path appears exactly once. 
For each $i \in I$, let $\mya_{i}(\mathcal{P}) := \rule[-1.75mm]{0.2mm}{5.5mm}\{j \in [1,k]_{\mathbb{Z}}\, |\, \xelt_{j-1} \myarrow{i_{j}} \xelt_{j}\}\rule[-1.75mm]{0.2mm}{5.5mm}$ (the number of color $i$ `ascents' in the path $\mathcal{P}$), and similarly set $\myd_{i}(\mathcal{P}) := \rule[-1.75mm]{0.2mm}{5.5mm}\{j \in [1,k]_{\mathbb{Z}}\, |\, \xelt_{j} \myarrow{i_{j}} \xelt_{j-1}\}\rule[-1.75mm]{0.2mm}{5.5mm}$ (the number of color $i$ `descents').

Now let $J$ be a subset of $I$. 
We let $\mathfrak{M}_{J}(\mathcal{P}) := \big\{j^{\mysmallera_{j}(\mathcal{P})+\mysmallerd_{j}(\mathcal{P})}\big\}_{j \in J}$ denote the multiset of colors from the set $J$ used in the path $\mathcal{P}$. 
So, for each $j \in J$, $\mathfrak{M}_{J}(\mathcal{P})$ contains the $\big(\mya_{j}(\mathcal{P})+\myd_{j}(\mathcal{P})\big)$-multisubset $\{j,j,\ldots,j\}$.  
Call our path $\mathcal{P}$ a $J${\em -path} if $\mya_{i}(\mathcal{P}) = 0 = \myd_{i}(\mathcal{P})$ for all $i \in I \setminus J$. 
The $J${\em -component} $\comp_{J}(\xelt)$ of an $\xelt$ in $R$ is the digraph consisting of all $\yelt$ that can be reached from $\xelt$ by some $J$-path together with the colored and directed edges of all such paths. 
If $\selt$ and $\telt$ are within the same $J$-component of $R$, then the $J${\em -distance} $\dist^{(J)}(\selt,\telt)$ between $\selt$ and $\telt$ is the minimum path length achieved when all $J$-paths from $\selt$ to $\telt$ in $R$ are considered; write $\dist^{(J)}(\selt,\telt)=\infty$ if $\selt$ and $\telt$ are not in the same $J$-component. 
When $\dist^{(J)}(\selt,\telt)<\infty$, any minimum-length-achieving $J$-path is {\em shortest}.  
Clearly any shortest $J$-path is simple. 
When $J=I$, we sometimes omit the prefix `$J$-' when speaking of paths, distance, etc. 

Say a given path $\mathcal{P}$ from $\selt$ to $\telt$ in $R$ is a {\em mountain path} if for some $j \in [0,k]_{\mathbb{Z}}$ and $\uelt := \xelt_{j}$ we have $\selt=\xelt_{0} \rightarrow \xelt_{1} \rightarrow \cdots \rightarrow \uelt \leftarrow \cdots \leftarrow \xelt_{k} = \telt$ (with edge colors suppressed), in which case we call $\uelt$ the {\em apex} of the mountain path. 
Similarly, call $\mathcal{P}$ a {\em valley path} if for some $j \in [0,k]_{\mathbb{Z}}$ and $\velt := \xelt_{j}$ we have $\selt=\xelt_{0} \leftarrow \xelt_{1} \leftarrow \cdots \leftarrow \velt \rightarrow \cdots \rightarrow \xelt_{k} = \telt$ (with edge colors suppressed), in which case we call $\velt$ the {\em nadir} of the valley path. 
A {\em nondecreasing} (resp., {\em nonincreasing}) path from $\selt$ to $\telt$ is a mountain (resp., valley) path whose apex (resp., nadir) is $\telt$. 

\begin{figure}[ht]
\begin{center}
\setlength{\unitlength}{0.4in} 
\begin{picture}(4,4.25)
\put(2,0){\circle*{0.1}} 
\put(1,1){\circle*{0.1}}
\put(2,1){\circle*{0.1}} 
\put(3,1){\circle*{0.1}} 
\put(0,2){\circle*{0.1}} 
\put(2,2){\circle*{0.1}} 
\put(4,2){\circle*{0.1}} 
\put(1,3){\circle*{0.1}} 
\put(2,3){\circle*{0.1}} 
\put(3,3){\circle*{0.1}} 
\put(2,4){\circle*{0.1}} 
\put(2,0){\line(-1,1){2}}
\put(2,0){\line(1,1){2}}
\put(0,2){\line(1,1){2}}
\put(4,2){\line(-1,1){2}}
\put(1,1){\line(1,1){2}}
\put(3,1){\line(-1,1){2}}
\put(2,0){\line(0,1){1}}
\put(2,4){\line(0,-1){1}}
\put(2,1){\line(-2,1){2}}
\put(2,1){\line(2,1){2}}
\put(2,3){\line(-2,-1){2}}
\put(2,3){\line(2,-1){2}}
\put(-4,-0.75){\small {\bf \TopoBalanceFigure}\ \  A connected and topographically balanced poset that is not a lattice.}
\end{picture}
\end{center}

\noindent
\vspace*{0.05cm}
\hspace*{4.85cm}$\overline{\hspace*{7cm}}$
\end{figure}

If $R$ is a poset with partial order `$\leq$', then let $\mbox{\sffamily Vertices}(R)$ be the set of elements of $R$ and $\mbox{\sffamily Edges}(R)$ the set of {\em covering relations} in $R$ wherein $\selt \rightarrow \telt$ if and only if $\selt < \telt$ and, for all $\relt \in R$, we have $\relt=\telt$ whenever $\selt < \relt \leq \telt$. 
The resulting (acyclic) digraph is the {\em covering digraph} of $R$. 
We identify $R$ with its covering digraph, so our poset $R$ can be edge-colored (resp.\ vertex-colored). 
In such an edge-colored poset $R$, we write $\xelt \myarrow{i} \yelt$ for vertices $\xelt$ and $\yelt$ and some color $i \in I$ if and only if $\xelt \rightarrow \yelt$ is a covering relation in $R$ and $\ecolor_{R}(\xelt \rightarrow \yelt) = i$. 
In figures, when such an edge is depicted without an arrowhead, the implied direction is `up'. 
Our poset $R$ is {\em topographically balanced} if (1) whenever  $\velt \rightarrow \selt$ and $\velt \rightarrow \telt$ for distinct $\selt$ and $\telt$ in $R$, 
then there exists a unique $\uelt$ in $R$ such that $\selt \rightarrow \uelt$ and $\telt \rightarrow \uelt$, and (2) whenever  $\selt \rightarrow \uelt$ and $\telt \rightarrow \uelt$ for distinct $\selt$ and $\telt$ in $R$, then there exists a unique $\velt$ in $R$ such that $\velt \rightarrow \selt$ and $\velt \rightarrow \telt$.  
See \TopoBalanceFigure\ for an example.  

A {\em rank function} on the poset $R$ is a surjective function $\rho : R \longrightarrow [0,l]_{\mathbb{Z}}$ (where $l \in [0,\infty)_{\mathbb{Z}}$) with the property that if $\selt \rightarrow \telt$ in $R$, then $\rho(\selt) + 1 = \rho(\telt)$. 
If such a rank function $\rho$ exists, then $R$ is a {\em ranked} poset and $l$ is the {\em length} of $R$ with respect to $\rho$. 
For any $\xelt \in R$ and $J \subseteq I$, let $\rho_{J}$ denote the rank function on a $J$-component $\comp_{J}(\xelt)$ of $R$, and let $l_{J}(\xelt)$ denote the length of $\comp_{J}(\xelt)$. 

We will make use of the following method by which any path $\mathcal{P}$ between vertices of a connected and topographically balanced poset $R$ can be converted to a mountain (respectively, valley) path $\mathcal{P}_{\mbox{\em \tiny mountain}}$ (resp., $\mathcal{P}_{\mbox{\em \tiny valley}}$). 
This is accomplished via our {\bf mountain{\_}path{\_}construction} algorithm. 
Suppose $\mathcal{P} = (\selt = \xelt_{0}, \xelt_{1}, \ldots, \xelt_{k-1}, \xelt_{k} = \telt)$ is a path between distinct elements $\selt$ and $\telt$ of $R$. 
We modify $\mathcal{P}$ as follows: 
\begin{enumerate}
\item If there is no $j \in [1,k-1]_{\mathbb{Z}}$ such that $\xelt_{j-1} \leftarrow \xelt_{j} \rightarrow \xelt_{j+1}$, then return $\mathcal{P}$ as $\mathcal{P}_{\mbox{\em \tiny mountain}}$. 
\item Otherwise, let $j \in [1,k-1]_{\mathbb{Z}}$ be least such that $\xelt_{j-1} \leftarrow \xelt_{j} \rightarrow \xelt_{j+1}$. 
If $\xelt_{j-1} = \xelt_{j+1}$, then form a new, and shorter, path $\mathcal{P}'$ from $\selt$ to $\telt$ by removing $\xelt_{j}$ and $\xelt_{j+1}$ from path $\mathcal{P}$, so that $\xelt_{j+2}$ now immediately succeeds $\xelt_{j-1}$ if $j+2 \leq k$; proceed now to Step 3. 

Assume now that $\xelt_{j-1} \not= \xelt_{j+1}$. 
Let $\xelt$ be the unique element from $R$ such that $\xelt_{j-1} \rightarrow \xelt \leftarrow \xelt_{j+1}$. 
We consider three cases, which are exhaustive and mutually exclusive: (a) $\xelt = \xelt_{k}$ for some $k \leq j-2$; (b) $\xelt = \xelt_{k}$ for some $k \geq j+2$; and (c) $\xelt$ is not visited by $\mathcal{P}$. 
In case (a), form a new path $\mathcal{P}'$ from $\selt$ to $\telt$ that is shorter than $\mathcal{P}$ by deleting $\xelt_{k+1}, \ldots, \xelt_{j}$ from $\mathcal{P}$, so $\xelt_{j+1}$ now immediately succeeds $\xelt_{k}$, and proceed to Step 3.  
Similarly, consider the case (b). 
In this case, form a new path $\mathcal{P}'$ from $\selt$ to $\telt$ that is shorter than $\mathcal{P}$ by deleting $\xelt_{j}, \ldots, \xelt_{k-1}$ from $\mathcal{P}$ so that $\xelt_{k}$ now immediately succeeds $\xelt_{j-1}$, and proceed to Step 3. 
Last, for case (c), set $\xelt_{j}' := \xelt$. 
Form a new path $\mathcal{P}'$ from $\mathcal{P}$ by replacing $\xelt_{j}$ with $\xelt_{j}'$.
Observe that $\mathcal{P}'$ is a path from $\selt$ to $\telt$ that is the same length as $\mathcal{P}$. 
Now go to Step 3.
\item Return to the first step of the process, using $\mathcal{P}'$ as $\mathcal{P}$. 
\end{enumerate} 
The `mountain-ization' $\mathcal{P}_{\mbox{\em \tiny mountain}}$ returned by {\bf mountain{\_}path{\_}construction} has the obvious properties that it is a a mountain path from $\selt$ to $\telt$ and is no longer than $\mathcal{P}$. 
The `valley-ization' $\mathcal{P}_{\mbox{\em \tiny valley}}$, which is a valley path from $\selt$ to $\telt$, can be obtained from $\mathcal{P}$ via a {\bf valley{\_}path{\_}construction} algorithm as follows: Apply {\bf mountain{\_}path{\_}construction} to the corresponding path $\mathcal{P}^{*}$ in $R^{*}$ to get $\mathcal{P\, }^{*}_{\mbox{\em $\!\!\!$\tiny mountain}}$, and then let $\mathcal{P}_{\mbox{\em \tiny valley}}$ be the path in $R$ corresponding to $\mathcal{P\, }^{*}_{\mbox{\em $\!\!\!$\tiny mountain}}$. 
Notice that if $\mathcal{P}$ is already a mountain (respectively, valley) path, then $\mathcal{P} = \mathcal{P}_{\mbox{\em \tiny mountain}}$ (resp., $\mathcal{P} = \mathcal{P}_{\mbox{\em \tiny valley}}$). 
Further, observe that if $\mathcal{P}$ is a shortest path from $\selt$ to $\telt$, then each of $\mathcal{P}_{\mbox{\em \tiny mountain}}$ and $\mathcal{P}_{\mbox{\em \tiny valley}}$ are also shortest.

The next result is a straightforward exercise and our first application of the preceding algorithm.

\noindent
{\bf \MountainValleyLemma}\ \ {\sl A connected and topographically balanced poset $R$ is ranked with a unique rank function we name $\rho$. 
Moreover, let $\mathcal{P}$ be a shortest path in $R$ from an element $\selt$ to another element $\telt$, let $\aelt$ be the apex of the mountain-ization} $\mathcal{P}_{\mbox{\em \tiny mountain}}$, {\sl and let $\nelt$ be the nadir of the valley-ization} $\mathcal{P}_{\mbox{\em \tiny valley}}$.  
{\sl Then} $\dist(\selt,\telt) = 2\rho(\aelt)-\rho(\selt)-\rho(\telt) = \rho(\selt)+\rho(\telt)-2\rho(\nelt)$. 

{\em Proof.} To demonstrate the claims of the last sentence of the result statement, we assume $R$ is ranked with unique rank function $\rho$. 
Then $\selt \leq \aelt$ and $\telt \leq \aelt$ $\Longrightarrow$ $\dist(\selt,\aelt) = \rho(\aelt)-\rho(\selt)$ and $\dist(\aelt,\telt) = \rho(\aelt)-\rho(\telt)$; since $\mathcal{P}$ is shortest, we have $\dist(\selt,\telt) = \dist(\selt,\aelt) + \dist(\aelt,\telt) = 2\rho(\aelt)-\rho(\selt)- \rho(\telt)$. 
Similarly see that $\dist(\selt,\telt) = \rho(\selt)+\rho(\telt) - 2\rho(\nelt)$. 

Now we demonstrate that $R$ is ranked. 
We begin by saying why $R$ must have a unique maximal element and a unique minimal element. 
Well, a maximal element exists since $R$ is finite. 
Now suppose there are two distinct, and therefore incomparable, maximal elements.  
Then any shortest path between them can be mountainized in order to produce an element that is above both maximal elements, contradicting their maximality. 
So, there exists a unique maximal element. 
Similarly see that there exists a unique minimal element, which, from here on, we call $\melt$. 

Next we show that for any $\xelt \in R$, any nondecreasing path from $\melt$ up to $\xelt$ has the same length as any shortest path from $\melt$ to $\xelt$. 
We argue this by induction on the distance of $\xelt$ from $\melt$. 
Of course, when $\dist(\melt,\xelt) = 0$, then $\xelt = \melt$, so any nondecreasing path from $\melt$ to $\xelt$ must have length zero, which is shortest. 
For our inductive hypothesis, suppose that for some nonnegative integer $r-1$, it is the case that whenever $\dist(\melt,\xelt) \leq r-1$, then any nondecreasing path from $\melt$ up to $\xelt$ has length equal to $\dist(\melt,\xelt)$. 
Now let $\mathcal{P} = (\melt = \xelt_{0},\ldots,\xelt_{p}=\xelt)$ be a nondecreasing path from $\melt$ up to some element $\xelt$, where $\dist(\melt,\xelt) = r$. 
In particular, $p \geq r$. 

Say $\mathcal{S} = (\melt = \selt_{0},\ldots,\selt_{r}=\xelt)$ is a shortest path from $\melt$ to $\xelt$. 
Then $\mathcal{S}_{\mbox{\em \tiny valley}}$ is a length $r$ valley path from $\melt$ to $\xelt$, and, since $\melt$ is the unique minimal element of $R$, the nadir of $\mathcal{S}_{\mbox{\em \tiny valley}}$ must also be $\melt$. 
In particular, $\mathcal{S}_{\mbox{\em \tiny valley}}$ is a length $r$ nondecreasing path from $\melt$ up to $\xelt$. 
Notice that for any $i \in [0,r]_{\mathbb{Z}}$, $(\melt = \selt_{0},\ldots,\selt_{i})$ is a shortest nondecreasing path from $\melt$ up to $\selt_{i}$. 
So, if $\xelt_{p-1} = \selt_{r-1}$, then our induction hypothesis applies to $\selt_{r-1}$, and therefore $p-1 = r-1$, i.e.\ $p=r$. 
That is, $\pathlength(\mathcal{P})=r$. 
Now suppose $\xelt_{p-1}$ and $\selt_{r-1}$ are distinct. 
Since $\xelt_{p-1} \rightarrow \xelt \leftarrow \selt_{r-1}$, then there exists a unique $\relt$ in $R$ such that $\xelt_{p-1} \leftarrow \relt \rightarrow \selt_{r-1}$. 
If we take any nondecreasing path $\mathcal{Q}$ from $\melt$ up to $\relt$, then we can append $\selt_{r-1}$ to get a nondecreasing path from $\melt$ up to $\selt_{r-1}$. 
The latter path has length $r-1$ by our induction hypothesis, and therefore the given nondecreasing path $\mathcal{Q}$ from $\melt$ up to $\relt$ has length $r-2$. 
So, if we append $\xelt_{p-1}$ to $\mathcal{Q}$, we get a path of length $r-1$ from $\melt$ up to $\xelt_{p-1}$. 
Again by our induction hypothesis, $\dist(\melt,\xelt_{p-1}) = r-1$, so we can conclude that $(\melt = \xelt_{0},\ldots,\xelt_{p-1})$ has length $r-1$, i.e.\ $p=r$. 
That is, $\mathcal{P}$ has length $r$. 
This completes the induction argument. 

Now define $\rho: R \longrightarrow [0,\infty)_{\mathbb{Z}}$ by $\rho(\xelt) := \dist(\melt,\xelt)$. 
We claim that if $\xelt \rightarrow \yelt$, then $\rho(\xelt)+1=\rho(\yelt)$, so that, with a suitable restriction on the target set $[0,\infty)_{\mathbb{Z}}$, $\rho$ is a rank function. 
So suppose $\xelt \rightarrow \yelt$. 
Any nondecreasing path $\mathcal{P}$ from $\melt$ up to $\xelt$ has length $\dist(\melt,\xelt)$, as argued above. 
If we append $\yelt$ to the path $\mathcal{P}$, we get a nondecreasing path of length $\dist(\melt,\xelt)+1$ from $\melt$ up to $\yelt$. 
Therefore, $\dist(\melt,\yelt) = \dist(\melt,\xelt)+1$, which means $\rho(\xelt)+1=\rho(\yelt)$. 
So, $\rho$ is a rank function. 
Connectedness of $R$ implies $\rho$ is unique.\hfill\QED

In a lattice $L$, $\xelt \vee \yelt$ denotes the join (least upper bound) of  elements $\xelt$ and $\yelt$, and $\xelt \wedge \yelt$ denotes their meet (greatest lower bound). 
We use $\mathsf{min}(L)$ (respectively, $\mathsf{max}(L)$) to denote its unique minimal (resp.\, maximal) element. 
A lattice $L$ is {\em distributive} if for any $\mathbf{r}$, $\mathbf{s}$, and $\mathbf{t}$ in $L$ it is the case that $\mathbf{r} \vee(\mathbf{s} \wedge \mathbf{t})=(\mathbf{r} \vee \mathbf{s}) \wedge (\mathbf{r} \vee \mathbf{t})$ and $\mathbf{r} \wedge (\mathbf{s} \vee \mathbf{t})=(\mathbf{r} \wedge \mathbf{s}) \vee (\mathbf{r} \wedge \mathbf{t})$.
A lattice is {\em modular} if it satisfies either of the equivalent conditions from the first sentence of the following proposition. 
This elementary result draws from precursor results in Theorem 1 of \cite{Alvarez}, in Section 2 of \cite{DuffusRival} (particularly Lemmas 2.5 and 2.6), in Section 3 of \cite{HG} (particularly Corollary 3.14), in Section 3.3 of \cite{StanleyText}, etc. 
This proposition applies to distributive lattices as well, since any distributive lattice is also modular. 

\noindent 
{\bf \ModularLatticeProp}\ \ 
{\sl A lattice $L$ (possibly edge-colored) is topographically balanced if and only if $L$ is ranked with unique rank function $\rho$ satisfying} 
\[2\rho(\selt \vee \telt) - \rho(\selt) - \rho(\telt) = \rho(\selt) + \rho(\telt) - 2\rho(\selt \wedge \telt)\] 
{\sl for all $\selt, \telt \in L$. 
Assume now that these equivalent conditions hold. Let $l$ be the length of $L$ with respect to $\rho$, and let $\mathcal{P}$ be a shortest path in $L$ from an element $\selt$ to an element $\telt$. Then} 
\[\displaystyle \mathsf{dist}(\selt,\telt) = \pathlength(\mathcal{P}) = 2\rho(\selt \vee \telt) - \rho(\selt) - \rho(\telt) = \rho(\selt) + \rho(\telt) - 2\rho(\selt \wedge \telt)<\infty.\] 
{\sl In particular, a mountain path from $\selt$ to $\telt$ is a shortest path if and only if its apex is $\selt \vee \telt$, and a valley path from $\selt$ to $\telt$ is shortest if and only if its nadir is $\selt \wedge \telt$. 
We have $\mathsf{dist}(\selt,\telt) \leq l$ for all $\selt$ and $\telt$ in $L$. 
Moreover, $\mathsf{dist}(\selt,\telt) = l$ if and only if $\selt \vee \telt = \mathsf{max}(L)$ and $\selt \wedge \telt = \mathsf{min}(L)$.} 

{\em Proof.} The first equivalence follows from Proposition 3.3.2 of \cite{StanleyText}. 
The last two sentences follow from Theorem 1 of \cite{Alvarez}. 
The formula for $\mathsf{dist}(\selt,\telt)$ can be obtained from Corollary 3.14 of \cite{HG} or Lemma 2.6 of \cite{DuffusRival}. 

Here is another proof of the above formula for $\mathsf{dist}(\selt,\telt)$ that also establishes the claim that immediately succeeds said formula in our proposition statement. 
Assume $L$ is modular. 
Suppose a path $\mathcal{P} = (\xelt_{0} = \selt, \xelt_{1}, \ldots, \xelt_{k-1}, \xelt_{k} = \telt)$ from $\selt$ to some distinct $\telt$ is shortest. 
Let $\aelt$ be the apex of the mountain-ization $\mathcal{P}_{\mbox{\em \tiny mountain}}$ of $\mathcal{P}$. 
Then $\pathlength(\mathcal{P}) = \dist(\selt,\telt) = 2\rho(\aelt)-\rho(\selt)-\rho(\telt)$, where the former equality is merely the definition of `shortest path' and the latter equality follows from \MountainValleyLemma. 
Now let $\mathcal{Q}$ be a path formed by a walk from $\selt$ up to $\selt \vee \telt$ followed by a walk from $\selt \vee \telt$ down to $\telt$. 
So, $\pathlength(\mathcal{Q}) = \left(\rho(\selt\vee \telt) - \rho(\selt)\right) + \left(\rho(\selt\vee \telt) - \rho(\telt)\right) = 2\rho(\selt \vee \telt) - \rho(\selt) - \rho(\telt)$. 
We see, then, that $2\rho(\aelt)-\rho(\selt)-\rho(\telt) = \dist(\selt,\telt) \leq \pathlength(\mathcal{Q}) = 2\rho(\selt \vee \telt) - \rho(\selt) - \rho(\telt)$, which means that $\rho(\aelt) \leq \rho(\selt \vee \telt)$. 
But since $\selt \leq \aelt$ and $\telt \leq \aelt$, then $\selt \vee \telt \leq \aelt$, and in particular $\rho(\selt \vee \telt) \leq \rho(\aelt)$. 
In any ranked poset, distinct elements with the same rank are incomparable. 
Therefore $(\selt \vee \telt) \leq \aelt$ together with $\rho(\selt \vee \telt) = \rho(\aelt)$ implies that $(\selt \vee \telt) = \aelt$. 
We conclude that $\pathlength(\mathcal{P}) = \dist(\selt,\telt) = 2\rho(\selt \vee \telt)-\rho(\selt)-\rho(\telt)$. 
Use the valley-ization $\mathcal{P}_{\mbox{\em \tiny valley}}$ in a similar way to conclude that $\pathlength(\mathcal{P}) = \dist(\selt,\telt) = \rho(\selt) + \rho(\telt) - 2\rho(\selt \wedge \telt)$.
As part of the preceding arguments, we observed that any mountain path from $\selt$ to $\telt$ with apex $\selt \vee \telt$ is a shortest path from $\selt$ to $\telt$, and we saw that any mountain path from $\selt$ to $\telt$ that is shortest must have $\selt \vee \telt$ as its apex. 
We can similarly establish the corresponding statements for valley paths from $\selt$ to $\telt$.\hfill\QED

When a modular or distributive lattice $L$ has its covering digraph edges colored by some set, we say this lattice is {\em diamond-colored} if, whenever  
\parbox{1.1cm}{\begin{center}
\setlength{\unitlength}{0.2cm}
\begin{picture}(4,3.5)
\put(2,0){\circle*{0.5}} \put(0,2){\circle*{0.5}}
\put(2,4){\circle*{0.5}} \put(4,2){\circle*{0.5}}
\put(0,2){\line(1,1){2}} \put(2,0){\line(-1,1){2}}
\put(4,2){\line(-1,1){2}} \put(2,0){\line(1,1){2}}
\put(0.75,0.55){\em \small k} \put(2.7,0.7){\em \small l}
\put(0.7,2.7){\em \small i} \put(2.75,2.55){\em \small j}
\end{picture} \end{center}} is an edge-colored subgraph of the covering digraph for $L$, then $i = l$ and $j = k$. 
We remark that in our research program studying poset models for semisimple Lie algebra representations and their associated Weyl group symmetric functions, diamond-colored modular and distributive lattices seem to arise naturally and seem to be (in some sense) highly efficient. 
Indeed, case-by-case investigations suggest at least the possibility that every irreducible simple Lie algebra representation and associated `Weyl bialternant' possesses at least one diamond-colored modular (though not necessarily distributive) lattice model. 


\begin{figure}[t]
{\Large 
$\mathbf{J}_{\mbox{\em \footnotesize color}}\left(\, \setlength{\unitlength}{0.6cm}
\begin{picture}(3,3)
\thicklines
\multiput(2,-1)(1,1){2}{\line(-1,1){2}}
\multiput(0,1)(1,-1){3}{\line(1,1){1}}
\put(0,1){\textcolor{red}{\circle*{0.3}}}
\put(-0.15,1.25){\em \scriptsize \textcolor{red}{1}}
\put(-0.15,0.45){\em \scriptsize \textcolor{red}{a}}
\put(1,2){\textcolor{blue}{\circle*{0.3}}}
\put(0.85,2.25){\em \scriptsize \textcolor{blue}{2}}
\put(0.85,1.45){\em \scriptsize \textcolor{blue}{b}}
\put(1,0){\textcolor{blue}{\circle*{0.3}}}
\put(0.85,0.25){\em \scriptsize \textcolor{blue}{2}}
\put(0.85,-0.55){\em \scriptsize \textcolor{blue}{c}}
\put(2,1){\textcolor{Green}{\circle*{0.3}}}
\put(1.85,1.25){\em \scriptsize \textcolor{Green}{3}}
\put(1.85,0.45){\em \scriptsize \textcolor{Green}{d}}
\put(2,-1){\textcolor{Green}{\circle*{0.3}}}
\put(1.85,-0.75){\em \scriptsize \textcolor{Green}{3}}
\put(1.85,-1.55){\em \scriptsize \textcolor{Green}{e}}
\put(3,0){\textcolor{orange}{\circle*{0.3}}}
\put(2.85,0.25){\em \scriptsize \textcolor{orange}{4}}
\put(2.85,-0.55){\em \scriptsize \textcolor{orange}{f}}
\end{picture}\, \right) = \hspace*{0.75cm} \rule[-3cm]{0mm}{6cm}
\setlength{\unitlength}{0.7cm}
\begin{picture}(3,0)
\thinlines
\put(1,-3.5){\textcolor{Green}{\line(0,1){1.5}}}
\put(0.9,-2.9){\textcolor{Green}{\em \scriptsize 3}}
\multiput(1,-2)(1,1){3}{\textcolor{orange}{\line(-1,1){1}}}
\multiput(0.4,-1.6)(1,1){3}{\textcolor{orange}{\em \scriptsize 4}}
\multiput(1,-2)(-1,1){2}{\textcolor{blue}{\line(1,1){1}}}
\multiput(1.4,-1.6)(-1,1){2}{\textcolor{blue}{\em \scriptsize 2}}
\multiput(2,-1)(-1,1){3}{\textcolor{red}{\line(1,1){1}}}
\multiput(2.4,-0.6)(-1,1){3}{\textcolor{red}{\em \scriptsize 1}}
\multiput(1,0)(1,1){2}{\textcolor{Green}{\line(-1,1){1}}}
\multiput(0.4,0.4)(1,1){2}{\textcolor{Green}{\em \scriptsize 3}}
\put(1,2){\textcolor{blue}{\line(0,1){1.5}}}
\put(0.9,2.6){\textcolor{blue}{\em \scriptsize 2}}
\put(1,-3.5){\circle*{0.2}}
\put(1.2,-3.65){\footnotesize $\emptyset$}
\put(1,-2){\circle*{0.2}}
\put(1,-2){\circle{0.35}}
\put(1.2,-2.15){\footnotesize $\langle \textcolor{Green}{e} \rangle$}
\put(0,-1){\circle*{0.2}}
\put(0,-1){\circle{0.35}}
\put(-0.85,-1.15){\footnotesize $\langle \textcolor{orange}{f} \rangle$}
\put(2,-1){\circle*{0.2}}
\put(2,-1){\circle{0.35}}
\put(2.2,-1.15){\footnotesize $\langle \textcolor{blue}{c} \rangle$}
\put(1,0){\circle*{0.2}}
\put(-0.35,-0.15){\footnotesize $\langle \textcolor{blue}{c},\!\textcolor{orange}{f} \rangle$}
\put(3,0){\circle*{0.2}}
\put(3,0){\circle{0.35}}
\put(3.2,-0.15){\footnotesize $\langle \textcolor{red}{a} \rangle$}
\put(0,1){\circle*{0.2}}
\put(0,1){\circle{0.35}}
\put(-0.85,0.85){\footnotesize $\langle \textcolor{Green}{d} \rangle$}
\put(2,1){\circle*{0.2}}
\put(2.2,0.85){\footnotesize $\langle \textcolor{red}{a},\!\textcolor{orange}{f} \rangle$}
\put(1,2){\circle*{0.2}}
\put(1.2,1.85){\footnotesize $\langle \textcolor{red}{a},\!\textcolor{Green}{d} \rangle$}
\put(1,3.5){\circle*{0.2}}
\put(1,3.5){\circle{0.35}}
\put(1.2,3.35){\footnotesize $\langle \textcolor{blue}{b} \rangle$}
\put(4.5,-0.25){\parbox{3.25in}{\footnotesize {\bf \JcolorFigure}\ \  The 10-element diamond-colored distributive lattice (`DCDL') depicted here is realized as the lattice of down-sets from a 6-element vertex-colored poset. 
The colored vertices of the latter poset are $\textcolor{red}{a}$ with color $\textcolor{red}{1}$; $\textcolor{blue}{b}$ and $\textcolor{blue}{c}$ with color $\textcolor{blue}{2}$; $\textcolor{Green}{d}$ and $\textcolor{Green}{e}$ with color $\textcolor{Green}{3}$; and $\textcolor{orange}{f}$ with color $\textcolor{orange}{4}$. 
Down-sets from this poset are notated as $\langle v_{1},\ldots,v_{k} \rangle$, where $v_{1},\ldots,v_{k}$ are the maximal vertices of the down-set.  
So, for example, $\langle \textcolor{red}{a},\!\textcolor{Green}{d} \rangle$ consists of the elements $\{\textcolor{red}{a},\!\textcolor{blue}{c},\!\textcolor{Green}{d},\!\textcolor{Green}{e},\!\textcolor{orange}{f}\}$.  
Since $\langle \textcolor{Green}{d} \rangle = \{\textcolor{blue}{c},\!\textcolor{Green}{d},\!\textcolor{Green}{e},\!\textcolor{orange}{f}\}$, then the DCDL edge $\langle \textcolor{Green}{d} \rangle {\textcolor{red}{\mylongarrow{1}}} \langle \textcolor{red}{a},\!\textcolor{Green}{d} \rangle$ has color {\textcolor{red}{1}} because $\langle \textcolor{red}{a},\!\textcolor{Green}{d} \rangle \setminus \langle \textcolor{Green}{d} \rangle = \{\textcolor{red}{a}\}$.
Join irreducible elements of the DCDL are circled.}} 
\end{picture}$
}
\end{figure}

Each diamond-coloring of a finite distributive lattice is naturally correlated with a vertex-coloring of its `compression poset' (aka poset of join irreducibles), in the following way. 
Suppose $L$ is a distributive lattice of `down-sets' (aka order ideals\footnote{A {\em down-set} or {\em order ideal} $\mathbf{i}$ from a poset $P$ is closed under `$\leq$' in the sense that $u \in \mathbf{i}$ whenever $u \leq v$ and $v \in \mathbf{i}$; for completeness, we remark that an {\em up-set} or {\em filter} $\mathbf{f}$ from $P$ is defined by the rule that $v \in \mathbf{f}$ whenever $u \leq v$ and $u \in \mathbf{f}$.}) taken from a vertex-colored poset $P$. 
Any covering relation (i.e.\ directed edge) $\xelt \rightarrow \yelt$ in the covering digraph of $L$ has color $i$ and is written $\xelt \myarrow{i} \yelt$  if and only if $\xelt \subset \yelt$ and the set difference $\yelt \setminus \xelt$ of down-sets $\xelt$ and $\yelt$ from $P$ is a singleton $\{v\}$ where vertex $v \in P$ is maximal element of $\yelt$ and has color $i$. 
In this case, we write $L = \mathbf{J}_{\mbox{\em \scriptsize color}}(P)$; it is easy to see that $L$ is diamond-colored and that the rank in $L$ of any down-set $\xelt$ from $P$ is just $|\xelt|$. 
On the other hand, starting with a diamond-colored distributive lattice $L$, assign color $i$ to a join irreducible\footnote{We say $\xelt$ from our distributive lattice $L$ is {\em join irreducible} if $\xelt$ covers exactly one element of $L$; for completeness, we remark that an element of $L$ is {\em meet irreducible} if it is covered by exactly one element of $L$.} element $\velt$ of $L$ if and only if $\uelt \myarrow{i} \velt$ when $\uelt$ is the unique element in $L$ covered by $\velt$. 
The induced-order subposet $P$ of $L$ consisting of join irreducible elements of $L$ is thusly vertex-colored, and we write $P = \mathbf{j}_{\mbox{\em \scriptsize color}}(L)$. 
What we call `The Fundamental Theorem of Finite Diamond-Colored Distributive Lattices (FTFDCDL)' is the following straightforward generalization of the `FTFDL' aka Birkhoff's Representation Theorem (for the latter, see Theorem 3.4.1 $\!\!/\!\!$ Proposition 3.4.2 of \cite{StanleyText} or Theorem 2.5 $\!\!/\!\!$ Corollary 2.6 of \cite{Aigner}): {\sl For any vertex-colored poset $P$ and diamond-colored distributive lattice $L$, we have} $P \cong \mathbf{j}_{\mbox{\em \scriptsize color}}(\mathbf{J}_{\mbox{\em \scriptsize color}}(P))$ {\sl and} $L \cong \mathbf{J}_{\mbox{\em \scriptsize color}}(\mathbf{j}_{\mbox{\em \scriptsize color}}(L))$. 
For an illustration of this result, see \JcolorFigure; for a formal proof see Theorem 2.8 of \cite{DonDistributive}.  

The next two results seem to be new, as there are no counterpart statements for monochromatic modular and distributive lattices. 
\ColorComponentTheorem\ seems to be interesting on its own (e.g.\ we easily obtain from it a generalization of a famous result due to R.\ P.\ Dilworth cf.\ \DilworthCorollary\ below), but mainly we require it for our development of \MoveMinGameTheorem.  

\noindent
{\bf \TechnicalLemma}\ \ {\sl Let $L$ be a diamond-colored modular lattice with edges colored by a set $I$, and take $\selt \leq \telt$ in $L$. 
Say $\mathcal{P} = (\selt = \xelt_{0} \myarrow{i_{1}} \xelt_{1} \myarrow{i_{2}} \cdots \myarrow{i_{p}} \xelt_{p}=\telt)$ and $\mathcal{Q} = (\selt = \yelt_{0} \myarrow{j_{1}} \yelt_{1} \myarrow{j_{2}} \cdots \myarrow{j_{q}} \yelt_{q}=\telt)$ are paths in $L$ from $\selt$ up to $\telt$, so in particular $\myd_{i}(\mathcal{P}) = 0 = \myd_{i}(\mathcal{Q})$ for all $i \in I$. 
Then $p=q$ and $\mathfrak{M}_{I}(\mathcal{P}) = \mathfrak{M}_{I}(\mathcal{Q})$. 
Moreover, if $\xelt_{1}$ and $\yelt_{p-1}$ are incomparable, then $i_{1} = j_{p}$.}

{\em Proof.} Since $L$ is ranked, then $p = q$. 
We use induction on the length $p$ of the given paths to prove both claims of the lemma statement.  
If $p = 0$, then there is nothing to prove.  
For our induction hypothesis, we assume the theorem statement holds whenever $p \leq m$ for some nonnegative integer $m$.  
Suppose now that $p = m+1$.  
We consider two cases: (1) $\xelt_{p-1} = \yelt_{p-1}$ and (2)  $\xelt_{p-1} \not= \yelt_{p-1}$. 
In case (1), the induction hypothesis applies to the paths $\selt = \xelt_{0} \myarrow{i_{1}} \xelt_{1} \myarrow{i_{2}} \xelt_{2} \myarrow{i_{3}} \cdots \myarrow{i_{p-1}} \xelt_{p-1} = \yelt_{p-1}$ and $\selt = \yelt_{0} \myarrow{j_{1}} \yelt_{1} \myarrow{j_{2}} \yelt_{2} \myarrow{j_{3}} \cdots \myarrow{j_{p-1}} \yelt_{p-1} = \xelt_{p-1}$.  
It follows that $\{i_{1}, i_{2}, \ldots, i_{p-1}\} \eqmulti \{j_{1}, j_{2}, \ldots, j_{p-1}\}$.  
Since in this case we have $i_{p} = j_{p}$, we conclude that $\{i_{1}, i_{2}, \ldots, i_{p-1}, i_{p}\} \eqmulti \{j_{1}, j_{2}, \ldots, j_{p-1}, j_{p}\}$.  
Also in this case, $\xelt_{1} \leq \xelt_{p-1} = \yelt_{p-1}$.  So $\xelt_{1}$ and $\yelt_{p-1}$ are comparable.  

In case (2), $\xelt_{p-1} \not= \yelt_{p-1}$. 
Let $\zelt := \xelt_{p-1} \wedge \yelt_{p-1}$.  
Since $\selt \leq \xelt_{p-1}$ and $\selt \leq \yelt_{p-1}$, it follows that $\selt \leq \zelt$.  
Consider a path $\selt = \welt_{0} \myarrow{k_{1}} \welt_{1} \myarrow{k_{2}} \welt_{2} \myarrow{k_{3}} \cdots \myarrow{k_{p-3}} \welt_{p-3} \myarrow{k_{p-2}} \welt_{p-2} = \zelt$. 
Since $L$ is a diamond-colored modular lattice, we have $\zelt \myarrow{j_{p}} \xelt_{p-1}$ and $\zelt \myarrow{i_{p}} \yelt_{p-1}$.  
Then by the induction hypothesis, we have $\{k_{1}, k_{2}, \ldots, k_{p-2}, j_{p}\} \eqmulti \{i_{1}, i_{2}, \ldots, i_{p-2}, i_{p-1}\}$ and $\{k_{1}, k_{2}, \ldots, k_{p-2}, i_{p}\} \eqmulti \{j_{1}, j_{2}, \ldots, j_{p-2}, j_{p-1}\}$.  
Then, \[\{i_{1}, i_{2}, \ldots, i_{p-2}, i_{p-1}, i_{p}\} \eqmulti \{k_{1}, k_{2}, \ldots, k_{p-2}, i_{p}, j_{p}\} \eqmulti \{j_{1}, j_{2}, \ldots, j_{p-2}, j_{p-1}, j_{p}\},\] as desired.  
Now suppose that $\xelt_{1}$ and $\yelt_{p-1}$ are incomparable.  
Suppose $\xelt_{1}$ and $\zelt$ are comparable.  
Then it must be the case that $\zelt < \xelt_{1}$.  
(Else, $\xelt_{1} \leq \zelt$ and $\zelt \leq \yelt_{p-1}$ means $\xelt_{1}$ and $\yelt_{p-1}$ are comparable.)  
Since $\xelt_{p-1} \not= \yelt_{p-1}$, then $p \geq 2$.  
If $p \geq 3$, then $\rho(\xelt_{1}) \leq \rho(\telt) - 2$, while $\rho(\zelt) = \rho(\telt) - 2$.  
This contradicts the fact that $\zelt < \xelt_{1}$, so when $\xelt_{1}$ and $\zelt$ are comparable, it must be the case that $p = 2$.  
If $p = 2$, then $\zelt = \selt$ and we have the diamond \parbox{2.5cm}{\begin{center}
\setlength{\unitlength}{0.2cm}
\begin{picture}(6.5,3.5)
\put(1,0){\begin{picture}(4,3.5)
\put(2,0){\circle*{0.5}} \put(0,2){\circle*{0.5}}
\put(2,4){\circle*{0.5}} \put(4,2){\circle*{0.5}}
\put(0,2){\line(1,1){2}} \put(2,0){\line(-1,1){2}}
\put(4,2){\line(-1,1){2}} \put(2,0){\line(1,1){2}}
\put(0.75,0.55){\em \tiny $i_{1}$} \put(2.7,0.7){\em \tiny $j_{1}$}
\put(0.7,2.7){\em \tiny $i_{p}$} \put(2.7,2.55){\em \tiny $j_{p}$}
\put(2.5,-0.75){\footnotesize $\selt$} \put(4.75,1.75){\footnotesize 
$\yelt_{p-1}$}
\put(2.5,4){\footnotesize $\telt$} \put(-3.5,1.75){\footnotesize $\xelt_{p-1}$}
\end{picture}} \end{picture} \end{center}} in $L$.  
From the diamond coloring property, we conclude that $i_{1} = j_{p}$.  
Suppose now that $\xelt_{1}$ and $\zelt$ are incomparable.  
Then we can apply the induction hypothesis to the paths $\selt = \welt_{0} \myarrow{k_{1}} \welt_{1} \myarrow{k_{2}} \welt_{2} \myarrow{k_{3}} \cdots \myarrow{k_{p-3}} \welt_{p-3} \myarrow{k_{p-2}} \welt_{p-2} = \zelt \myarrow{j_{p}} \xelt_{p-1}$ and $\selt = \xelt_{0} \myarrow{i_{1}} \xelt_{1} \myarrow{i_{2}} \xelt_{2} \myarrow{i_{3}} \cdots \myarrow{i_{p-1}} \xelt_{p-1}$.  
From this, we see that $i_{1} = j_{p}$.  
This completes the induction step, and the proof.\hfill\QED  

The preceding result is applied in our proof of the claim, stated as part of \ColorComponentTheorem\ below, that any $J$-component of a diamond-colored modular (distributive) lattice $L$ is itself a modular (distributive) lattice in a natural way. 
In the case that $L$ is distributive, the next theorem also provides a way to identify the compression posets of the $J$-components of $L$ as certain induced-order subposets of its vertex-colored compression poset.  
To set up this latter statement, we need the following notions. 

Suppose $P$ is a poset with vertex coloring function $\vcolor: P \longrightarrow I$. 
Let $Q$ be a subposet of $P$ in the induced order and with the inherited vertex coloring.  
Call $Q$ a $J$-{\em subordinate} of $P$ if ({\em i}) $\vcolor(Q) \subseteq J$; ({\em ii}) there is a down-set $\relt$ from $P$ such that $\relt \cap Q = \emptyset$ and $\relt \cup Q$ is a down-set from $P$, i.e.\ $Q = \selt \setminus \relt$ where for some down-set $\selt$ containing the down-set $\relt$; and ({\em iii}) $\vcolor(v) \in I \setminus J$ whenever $v$ is a maximal (respectively, minimal) element of $\relt$ (resp.\ $P \setminus (\relt \cup Q)$). 
Regarding $\xelt$ to be a down-set from $P$, suppose $D_{1}$ and $D_{2}$ are any two subsets of color $J$ vertices in $\xelt$ such that $\xelt \setminus D_{1}$ and $\xelt \setminus D_{2}$ are both down-sets from $P$. 
It is easy to see that $\xelt \setminus (D_{1} \cup D_{2})$ is also a down-set from $P$. 
So, we may define $D_{J}(\xelt)$ be the uniquely largest subset of color $J$ vertices in $\xelt$ such that $\xelt \setminus D_{J}(\xelt)$ is a down-set from $P$.  
(Note that `$D$' is short for `delete'.) 
Similarly let $A_{J}(\xelt)$ be the uniquely largest subset of color $J$ vertices in the up-set $P \setminus \xelt$ such that $\xelt \cup A_{J}(\xelt)$ is an up-set from $P$.  
(`$A$' here is short for `add'.) 
Let $Q_{J}(\xelt) := A_{J}(\xelt) \cup D_{J}(\xelt)$. 
View $Q_{J}(\xelt)$ as a subposet in the induced order and with the inherited vertex coloring.

\noindent
{\bf \ColorComponentTheorem}\ \ {\sl Let $L$ be a diamond-colored modular lattice with edges colored by a set $I$. 
Let $J \subseteq I$ and $\xelt \in L$. 
(1) Then} $\comp_{J}(\xelt)$ {\sl is the covering digraph for a diamond-colored modular lattice and is a sublattice of $L$ in the sense that $\selt \vee_{L} \telt$ and $\selt \wedge_{L} \telt$ are in} $\comp_{J}(\xelt)$ {\sl whenever} $\selt, \telt \in \comp_{J}(\xelt)$. 
{\sl Moreover}, $\comp_{J}(\xelt)$ {\sl is distributive if $L$ is distributive.} 
{\sl (2) Now suppose $L$ is distributive with} $P := \mathbf{j}_{\mbox{\em \scriptsize color}}(L)$. 
{\sl Then $Q_{J}(\xelt)$ is a $J$-subordinate of $P$ with} $\comp_{J}(\xelt) \cong \mathbf{J}_{\mbox{\em \scriptsize color}}(Q_{J}(\xelt))$ {\sl and equivalently} $\mathbf{j}_{\mbox{\em \scriptsize color}}(\comp_{J}(\xelt)) \cong Q_{J}(\xelt)$.  
{\sl Moreover, each $J$-subordinate of $P$ is precisely $Q_{J}(\yelt)$ for some $\yelt \in L$.} 

{\em Proof.} The last sentence of part {\sl (1)} will follow immediately once we establish all prior claims in part {\sl (1)} of the theorem statement. 
Let $K$ be the edge-colored digraph $\comp_{J}(\xelt)$, and define on $K$ the relation $\leq_{K}$ wherein $\selt \leq_{K} \telt$ if and only if there is a path $(\selt=:\yelt_{0} \myarrow{j_{1}} \yelt_{1} \myarrow{j_{2}} \cdots \myarrow{j_{p}} \yelt_{p}:=\telt)$ such that each $\yelt_{i}$ is in $K$ and each $j_{i}$ is in $J$. 
Of course, $\leq_{K}$ is just the partial order on $K$ induced by $L$, and the edges $\selt \myarrow{j} \telt$ in $K$ are exactly the covering relations. 

Next, for arbitrary $\selt$ and $\telt$ in $K$, we prove, by inducting on the distance $\dist^{(J)}(\selt,\telt)$ between vertices of $K$, that any shortest path from $\selt$ to $\telt$ in $K$ is also a shortest path in $L$ and that $K$ contains both $\selt \vee_{L} \telt$ and $\selt \wedge_{L} \telt$. 
From this is will follow that the sublattice $K$ is the covering digraph for a diamond-colored modular lattice. 
When $\dist^{(J)}(\selt,\telt) = 0$, then $\selt=\telt$. 
The only shortest path in $K$ from $\selt$ to $\telt$ is the single-node path $(\selt)$, which is clearly shortest in $L$; also, $\selt \vee_{L} \telt = \selt = \selt \wedge_{L} \telt$. 
For our inductive hypothesis, we suppose that $p-1$ is a nonnegative integer such that for all $\selt', \telt' \in K$ with $\dist^{(J)}(\selt',\telt') \leq p-1$, then any shortest path from $\selt'$ to $\telt'$ in $K$ is also a shortest path in $L$ and that $K$ contains both $\selt' \vee_{L} \telt'$ and $\selt' \wedge_{L} \telt'$. 

Suppose now that $\dist^{(J)}(\selt,\telt) = p$. 
Let $\mathcal{P} := (\selt=\yelt_{0},\ldots,\yelt_{p}=\telt)$ be a shortest path in $K$ from $\selt$ to $\telt$. 
Assume for the moment that $\yelt_{p-1} \mylongbackarrow{j_{p}} \yelt_{p}$, and let $\mathcal{P}' := (\selt=\yelt_{0},\ldots,\yelt_{p-2},\yelt_{p-1})$. 
Then $\mathcal{P}'$ from $\selt$ to $\yelt_{p-1}$ has length $p-1$ and must be a shortest path in $K$. 
By our inductive hypothesis, $\mathcal{P'}$ is also a shortest path in $L$. 
Since $L$ is topographically balanced and $\mathcal{P}'$ is a shortest path in $L$ from $\selt$ to $\yelt_{p-1}$, we can apply our {\bf mountain{\_}path{\_}construction} algorithm to $\mathcal{P}'$ to obtain a shortest mountain path $\mathcal{P}'_{\mbox{\em \tiny mountain}}$ from $\selt$ to $\yelt_{p-1}$ whose apex $\aelt = \selt \vee_{L} \yelt_{p-1}$. 
At each application of Step 2 during the execution of the {\bf mountain{\_}path{\_}construction} algorithm, diamond-coloring of $L$ ensures that the multiset of path colors $\{j_{1},\ldots,j_{p-1}\}$ is unchanged. 
So, the multiset of edge colors $\{j_{1},\ldots,j_{p-1}\}$ for $\mathcal{P}'$ is the same as the multiset of edge colors for $\mathcal{P}'_{\mbox{\em \tiny mountain}}$. 
Thus, $\mathcal{P}'_{\mbox{\em \tiny mountain}}$ is a path in $K$, and in particular $\aelt \in K$. 
We let $\aelt' := \selt \vee_{L} \telt$ and observe that $\aelt' \leq \aelt$ since $\selt \leq \aelt$ and $\telt < \yelt_{p-1} \leq \aelt$. 
Take any path $\mathcal{Q}$ in $L$ from $\selt$ up to $\aelt$ that goes through $\aelt'$. 
Another path $\mathcal{Q}'$ that starts at $\selt$ and ascends to $\aelt$ is the initial part of $\mathcal{P}'_{\mbox{\em \tiny mountain}}$. 
Since $\mathcal{Q}'$ is a path in $K$, its edge colors are from the set $J$. 
By \TechnicalLemma, the multiset of edge colors in the path $\mathcal{Q}$ coincides with the multiset of edge colors in the path $\mathcal{Q}'$. 
Therefore, $\mathcal{Q}$ is a path in $K$, and so $\aelt' \in K$. 
Let $\mathcal{R}'$ be the path that starts at $\aelt$ and descends to $\yelt_{p-1}$ on edges from $\mathcal{P}'_{\mbox{\em \tiny mountain}}$ and then descends from $\yelt_{p-1}$ to $\telt$ along our given edge of color $j_{p}$. 
Further, let $\mathcal{R}$ be any path in $L$ from $\aelt$ down to $\telt$ that goes through $\aelt'$. 
Again by \TechnicalLemma, the multiset of edge colors in the path $\mathcal{R}$ coincides with the multiset of edge colors in the path $\mathcal{R}'$. 
Therefore $\mathcal{R}$ is a path in $K$. 
Consider the concatenation of $\mathcal{Q}$ and $\mathcal{R}$ but with any edges in the interval $[\aelt',\aelt]_{K}$ removed; the result is a mountain path $\mathcal{M}$ from $\selt$ to $\telt$ in $K$ whose apex is $\aelt' = \selt \vee_{L} \telt$. 
Moreover, this mountain path $\mathcal{M}$ is shortest in $L$ by \ModularLatticeProp, so $\dist_{L}(\selt,\telt)$ is the length of $\mathcal{M}$. 
So, $\mathcal{M}$ must also be shortest in $K$, and therefore $\mathcal{M}$ has length $p = \dist^{(J)}(\selt,\telt)$. 
That is, the path $\mathcal{P}$ of length $p$ from $\selt$ to $\telt$ in $K$ that we were originally given is a shortest path in $L$. 
Moreover, if we apply our {\bf valley{\_}path{\_}construction} algorithm to $\mathcal{M}$, the valley path $\mathcal{M}_{\mbox{\em \tiny valley}}$ we obtain is a shortest path in $K$ (since, due to earlier observations concerning Step 2 of the algorithm and diamond-coloring, its multiset of edge colors is the same as the multiset of edge colors for $\mathcal{M}$) and is shortest in $L$ (a feature of the algorithm is that it preserves shortest path lengths). 
Therefore by \ModularLatticeProp, its nadir is $\selt \wedge_{L} \telt$, which is therefore an element of $K$.

To summarize our work so far on the inductive step of our proof of part {\sl (1)}, we have shown that if $\mathcal{P} := (\selt=\yelt_{0},\ldots,\yelt_{p}=\telt)$ is a shortest path in $K$ from $\selt$ to $\telt$ such that $\yelt_{p-1} \mylongbackarrow{j_{p}} \yelt_{p}=\telt$, then we obtain that $\mathcal{P}$ a shortest path in $L$ from $\selt$ to $\telt$, that $\selt \vee_{L} \telt \in K$, and that $\selt \wedge_{L} \telt \in K$. 
An entirely similar argument can be given in the case that, for our given path $\mathcal{P}$, we have $\yelt_{p-1} \mylongarrow{j_{p}} \yelt_{p}=\telt$. 
This completes the inductive step, and the proof of part {\sl (1)}. 

For {\sl (2)}, let $Q = Q_{J}(\xelt)$.  
With $\relt := \xelt \setminus D_{J}(\xelt) = \mathsf{min}(\comp_{J}(\xelt))$ and $\relt \cup Q = \mathsf{max}(\comp_{J}(\xelt))$, it is easy to see that $Q$ meets the criteria for a $J$-subordinate of $P$. 
Now let $\phi: \mathbf{J}_{\mbox{\em \scriptsize color}}(Q_{J}(\xelt)) \longrightarrow \comp_{J}(\xelt)$ be given by $\phi(\yelt) = \yelt \cup \xelt$. 
It is routine to check that $\phi$ is an edge and edge-color preserving mapping between diamond-colored distributive lattices.
From $\comp_{J}(\xelt) \cong \mathbf{J}_{\mbox{\em \scriptsize color}}(Q_{J}(\xelt))$, obtain the congruence $\mathbf{j}_{\mbox{\em \scriptsize color}}(\comp_{J}(\xelt)) \cong Q_{J}(\xelt)$ by applying the Fundamental Theorem for Finite Diamond-colored Distributive Lattices. 
Suppose now we are given some $J$-subordinate $Q$ of $P$. 
Let $\yelt$ be the down-set $\relt$ of part ({\em ii}) of the definition of $J$-subordinate.  
Then it is easy to see that $Q = Q_{J}(\yelt)$.\hfill\QED

We pause briefly here to note an interesting order-theoretic consequence of the preceding theorem. 
A famous result of R.\ P.\ Dilworth \cite{Dilworth} states that in any modular lattice and for any positive integer $k$, the number of elements which cover exactly $k$ vertices is the same as the number of elements covered by exactly $k$ vertices (for some discussion, see for example the solution to Ch.\ 3 Exercise 101 (d) in \cite{StanleyText}). 
In particular, the number of join irreducible elements is the same as the number of meet irreducible elements. 
Let $J$ be a subset of our set of colors $I$. 
If $j \in J$ and $\xelt \myarrow{j} \yelt$ is an edge within some diamond-colored modular lattice $L$, then say that $\xelt$ is a $J${\em -descendant} of $\yelt$ and that $\yelt$ is a $J${\em -ascendant} of $\xelt$. 
Let $\#_{\mbox{\tiny $J$-asc}}(\xelt)$ denote the number of $J$-ascendants of $\xelt$, and let $\#_{\mbox{\tiny $J$-desc}}(\yelt)$ denote the number of $J$-descendants of $\yelt$. 

\noindent 
{\bf \DilworthCorollary}\ \ {\sl Let $L$ be a diamond-colored modular lattice, let $J \subseteq I$, and let $k \in [1,\infty)_{\mathbb{Z}}$. 
Then} $\rule[-1.75mm]{0.2mm}{5.5mm}\{\xelt \in L\, |\, \#_{\mbox{\tiny $J$-asc}}(\xelt) = k\}\rule[-1.75mm]{0.2mm}{5.5mm} = \rule[-1.75mm]{0.2mm}{5.5mm}\{\yelt \in L\, |\, \#_{\mbox{\tiny $J$-desc}}(\yelt) = k\}\rule[-1.75mm]{0.2mm}{5.5mm}$.

{\em Proof.} By \ColorComponentTheorem, Dilworth's result generalizes to our multi-colored setting.\hfill\QED 

We close this section by informally describing another application of \ColorComponentTheorem\ related to diamond-colored modular and/or distributive lattice models of finite-dimensional complex semisimple Lie algebra representations (and their companion Weyl group symmetric functions). 
See \cite{DonSupp} for an orientation to the language of this paragraph. 
Say $I$ is the indexing set for Chevalley generators of a finite-dimensional complex semisimple Lie algebra we momentarily denote by $\mathfrak{g}(I)$, and suppose $L$ is a DCML (DCDL) supporting graph for a representation of $\mathfrak{g}(I)$. 
For any subset $J \subseteq I$, denote by $\mathfrak{g}(J)$ the Levi subalgebra of $\mathfrak{g}(I)$ whose generators are those indexed by the subset $J$. 
Then any $J$-component $K$ of $L$ is, naturally, a DCML (DCDL) supporting graph for a representation of $\mathfrak{g}(J)$. 
In the distributive lattice case, the compression poset $\mathbf{j}_{\mbox{\em \scriptsize color}}(K)$ is therefore a $J$-subordinate of $\mathbf{j}_{\mbox{\em \scriptsize color}}(L)$.

\vspace{1ex} 
\noindent 
{\Large \bf \MoveMinGame\hspace*{0.15in}Lattice-theoretic solutions to certain move-minimizing puzzles}

\noindent 
\renewcommand{\thefootnote}{5} 
As observed in Section~\Intro, a combinatorial puzzle of minimizing moves between objects coincides with analyzing shortest paths between nodes in some kind of edge-colored directed graph. 
We make this more precise as follows. 
Let $\mathcal{G}$ be a finite simple directed graph with edges colored by some finite index set $I$. 
The vertices of $\mathcal{G}$ are to be regarded as \underline{objects} while the directed and colored edges of $\mathcal{G}$ are to be considered \underline{moves}. (We use the adjective `legal', as in `legal move', to emphasize that a given move is allowed within the context of the given puzzle.) 
The challenges posed in the sets of questions\footnote{Of course, someone analyzing a move-minimizing puzzle digraph might wish to add other move-minimizing questions to the list, such as ``What is the total number of shortest paths between any two given objects?''} (1a-b) and (2a-b) of Section \Intro\ constitute the {\em move-minimizing puzzle on} $\mathcal{G}$, and we call $\mathcal{G}$ itself the {\em puzzle digraph of objects and moves} ({\em puzzle digraph} for short) associated with the move-minimizing puzzle. 
Notice that if we are not interested in using colors to distinguish between moves, then we may simply look at the uncolored directed graph, in which case questions (1a-b) are the only relevant questions. 
We understand that the meaning of the term `explicit' in questions (1a-b) and (2a-b) is open to interpretation and would depend upon (for example) how explicit is the combinatorial description of the objects and moves of a given puzzle. 
The following theorem connects move-minimizing puzzles and diamond-colored modular and distributive lattices and indicates the level of explicitness we aim to achieve in answers to the move-minimizing questions we have posed. 
While there are aspects of this theorem that following directly from standard ideas, some aspects require non-standard ideas developed in Section \Background. 

\noindent 
{\bf \MoveMinGameTheorem}\ \ 
{\sl Suppose a move-minimizing puzzle digraph $\mathcal{G}$ with edges colored by a set $I$ is the covering digraph for some diamond-colored modular lattice $L$ whose rank function is $\rho$ and whose length is $l$. 
Let $\selt, \telt \in \mathcal{G}$. Then:}
\begin{enumerate}
\item[{\sl (1)}] {\sl (a) There exists an $I$-path in $\mathcal{G}$ from $\selt$ to $\telt$. 
More specifically, the distance from $\selt$ to $\telt$ in $\mathcal{G}$ is $\mathsf{dist}^{(I)}(\selt,\telt) = 2\rho(\selt \vee \telt) - \rho(\selt) - \rho(\telt) = \rho(\selt) + \rho(\telt) - 2\rho(\selt \wedge \telt)<\infty$. 
Indeed, $\mathsf{dist}^{(I)}(\selt,\telt) \leq l$, and, moreover, $\selt \wedge \telt = \mathsf{min}(L)$ and $\selt \vee \telt = \mathsf{max}(L)$ if and only if $\mathsf{dist}^{(I)}(\selt,\telt) = l$.}

 {\sl (b) If $L$ is distributive, let} $P := \mathbf{j}_{\mbox{\em \scriptsize color}}(L)$ {\sl be the affiliated vertex-colored compression poset and identify $\selt$ and $\telt$ as down-sets from $P$, so in particular $\selt$ and $\telt$ are subsets of $P$.  
We have the following expressions for certain quantities from part (1a): $\rho(\selt) = |\selt|$, $\rho(\telt) = |\telt|$, $\rho(\selt \vee \telt) = |\selt \cup \telt|$, and $\rho(\selt \wedge \telt) = |\selt \cap \telt|$.} 
{\sl Now let $\myl_{1} := |\selt \cup \telt| - |\selt|$ and $\myl_{2} := |\selt \cup \telt| - |\telt|$. 
Let $\xelt_{0} := \selt$. 
For $j \in [1,\myl_{1}]_{\mathbb{Z}}$, let $u$ be any minimal element of the set $(\selt \cup \telt) \setminus \xelt_{j-1}$ and declare $\xelt_{j} := \xelt_{j-1} \cup \{u\}$.  Let $\xelt_{\ell_{1}+\ell_{2}} := \telt$.  For $j \in [1,\myl_{2}]_{\mathbb{Z}}$, let $u$ be any minimal element of the set $(\selt \cup \telt) \setminus \xelt_{\ell_{1} + j}$ and declare $\xelt_{\ell_{1} + j - 1} := \xelt_{\ell_{1} + j} \cup \{u\}$.  
Then for $j \in [0,\myl_{1}+\myl_{2}]_{\mathbb{Z}}$, $\xelt_{j}$ is an order ideal from $P$, and $(\xelt_{0}, \xelt_{1}, \ldots , \xelt_{\ell_{1}+\ell_{2}})$ is a shortest path from $\selt$ to $\telt$ with $\xelt_{\ell_{1}} = \selt \cup \telt$. 
Similarly obtain a shortest path $(\selt = \yelt_{0}, \yelt_{1}, \ldots , \yelt_{k_{1}+k_{2}}=\telt)$ from $\selt$ to $\telt$ with $\yelt_{k_{1}} = \selt \cap \telt$, where $k_{1} = |\selt| - |\selt \cap \telt|$ and $k_{2} = |\telt| - |\selt \cap \telt|$.} 

\item[{\sl (2)}] {\sl (a) Say $J \subseteq I$.  
Let $\mathcal{P}$ be any shortest $I$-path from $\selt$ and $\telt$, and let $\mathcal{Q}$ be any path from $\selt$ to $\telt$ such that the total number of edges in $\mathcal{Q}$ with colors from the set $J$ is as small as possible, i.e.\ $\rule[-1.75mm]{0.2mm}{5.5mm}\, \mathfrak{M}_{J}(\mathcal{Q})\rule[-1.75mm]{0.2mm}{5.5mm}$ is minimized.  
Then $\mathfrak{M}_{J}(\mathcal{P}) = \mathfrak{M}_{J}(\mathcal{Q})$.}
{\sl Moreover, there exists a $J$-path in $\mathcal{G}$ from $\selt$ to $\telt$ if and only if} $\comp_{J}(\selt) = \comp_{J}(\telt)$; {\sl in this case, any shortest $I$-path from $\selt$ to $\telt$ is necessarily a $J$-path, and $\mathsf{dist}^{(J)}(\selt,\telt) = \mathsf{dist}^{(I)}(\selt,\telt)$.}

{\sl (b) Suppose that $L$ is distributive, as in part (1b).  Then} $\comp_{J}(\selt) = \comp_{J}(\telt)$ {\sl if and only if $Q_{J}(\selt) = Q_{J}(\telt)$ (i.e.\ we have equality of the corresponding $J$-subordinates of $P$).}  
{\sl Under these equivalent conditions, let $K$ be the edge-colored digraph} $\comp_{J}(\selt) = \comp_{J}(\telt)$ {\sl and a subposet of $L$ in the induced order.}  
{\sl Then the  $J$-subordinates $Q_{J}(\selt)$ and $Q_{J}(\telt)$ of $P$ both coincide with the set $Q:=(\selt \cup \telt)\setminus(\selt \cap \telt)$, and moreover} $K \cong \mathbf{J}_{\mbox{\em \scriptsize color}}(Q)$. 
{\sl Also, $\mathsf{dist}^{(I)}(\selt,\telt) = \mathsf{dist}^{(J)}(\selt,\telt) = 2\rho_{J}(\selt \vee \telt) - \rho_{J}(\selt) - \rho_{J}(\telt) = \rho_{J}(\selt) + \rho_{J}(\telt) - 2\rho_{J}(\selt \wedge \telt)<\infty$, where $\rho_{J}(\selt) = |\selt\setminus(\selt \cap \telt)|$, $\rho_{J}(\telt) = |\telt\setminus(\selt \cap \telt)|$, $\rho_{J}(\selt \vee \telt) = |(\selt \cup \telt)\setminus(\selt \cap \telt)|=l_{J}(\selt)$, and $\rho_{J}(\selt \wedge \telt) = 0$.}
 {\sl Further, one can locate a shortest $J$-path from $\selt$ to $\telt$ -- which is also a shortest $I$-path from $\selt$ to $\telt$ -- by applying the process specified in part (1b) above to the $J$-component $K$ and its compression poset $Q$.}
\end{enumerate}


{\em Proof.} Part {\sl (1a)} follows from \ModularLatticeProp. 
For {\sl (1b)}, it follows from the definitions that each set $\xelt_{j}$ is an order ideal and that $\xelt_{j-1}$ and $\xelt_{j}$ are neighbors in the covering digraph for $L$. 
As for existence of such a set $\xelt_{j}$, consider first the case that $j \in [1,\ell_{1}]_{\mathbb{Z}}$. 
Since $\selt \subseteq \selt \cup \telt$ and $\rho(\selt \cup \telt) - \rho(\selt) = |\selt \cup \telt|-|\selt|$, then the length $\ell_{1}$ of any path in $L$ from $\selt$ up to $\selt \cup \telt$ is $|\selt \cup \telt|-|\selt|$. 
Similarly, $\ell_{2} = |\selt \cup \telt|-|\telt|$. 
Since the path $(\xelt_{0}, \xelt_{1}, \ldots , \xelt_{\ell_{1}+\ell_{2}})$ from $\selt$ to $\telt$ has length $\myl_{1}+\myl_{2} = |\selt \cup \telt| - |\selt| + |\selt \cup \telt| - |\telt| = 2\rho(\selt \vee \telt) - \rho(\selt) - \rho(\telt)$, then, by \ModularLatticeProp, the path is shortest. 
Similar reasoning shows $(\selt = \yelt_{0}, \yelt_{1}, \ldots , \yelt_{k_{1}+k_{2}}=\telt)$ is a shortest path in $L$ from $\selt$ to $\telt$.

The claims of {\sl (2b)} follow from the definitions together with \ColorComponentTheorem\ and part {\sl (1)} of the theorem statement.  
For the claim in {\sl (2a)} that there exists a $J$-path in $\mathcal{G}$ from $\selt$ to $\telt$ if and only if $\comp_{J}(\selt) = \comp_{J}(\telt)$, simply observe that the binary relation on the vertices of $\mathcal{G}$ defined by whether two vertices are connected by a $J$-path is an equivalence relation. 
In this case, then by \ColorComponentTheorem, the $J$-component $K := \comp_{J}(\selt)=\comp_{J}(\telt)$ is a sublattice of $L$ with $\selt \vee \telt$ and $\selt \wedge \telt$ both residing in this $J$-component. 
In particular, there exists a $J$-path $\mathcal{R}$ from $\selt$ to $\telt$ that is also a shortest $I$-path, namely any path in $K$ that ascends from $\selt$ to $\selt \vee \telt$ and from there descends in $K$ to $\telt$. 
Now, any shortest $I$-path $\mathcal{S}$ from $\selt$ to $\telt$ can be modified to become a mountain path $\mathcal{S}_{\mbox{\em \tiny mountain}}$ through $\selt \vee \telt$, by our {\bf mountain{\_}path{\_}construction} algorithm, and it follows from the application of Step 2 of that algorithm that $\mathfrak{M}_{I}(\mathcal{S}) = \mathfrak{M}_{I}(\mathcal{S}_{\mbox{\em \tiny mountain}})$. 
Of course, by applying \TechnicalLemma\ to the `ascending' parts of $\mathcal{R}$ and $\mathcal{S}_{\mbox{\em \tiny mountain}}$ and then again to the `descending' part of these two mountain paths, we see that $\mathfrak{M}_{I}(\mathcal{R}) = \mathfrak{M}_{I}(\mathcal{S}_{\mbox{\em \tiny mountain}})$. 
Since $\mathcal{R}$ is a $J$-path, $\mathcal{S}_{\mbox{\em \tiny mountain}}$, and hence $\mathcal{S}$, must be $J$-paths. 
From this, we also see that $\mathsf{dist}^{(J)}(\selt,\telt) = \pathlength(\mathcal{R}) = \pathlength(\mathcal{S}) = \mathsf{dist}^{(I)}(\selt,\telt)$. 

Last, we demonstrate the claim, stated at the beginning of part {\sl (2a)}, that $\mathfrak{M}_{J}(\mathcal{P}) = \mathfrak{M}_{J}(\mathcal{Q})$. 
Here we make no assumption that $\selt$ and $\telt$ are connected by a $J$-path. 
As in the previous paragraph, the mountain path $\mathcal{P}_{\mbox{\em \tiny mountain}} = (\selt =: \xelt_{0},\ldots,\xelt_{p}:=\telt)$ is also a shortest path from $\selt$ to $\telt$ and $\mathfrak{M}_{I}(\mathcal{P}_{\mbox{\em \tiny mountain}}) = \mathfrak{M}_{I}(\mathcal{P})$. 
Then the apex of $\mathcal{P}_{\mbox{\em \tiny mountain}}$ is $\selt \vee \telt$. 
The mountain path $\mathcal{Q}_{\mbox{\em \tiny mountain}}$ need not be shortest, since $\mathcal{Q}$ is not necessarily shortest. 
Since Step 2 of the {\bf mountain{\_}path{\_}construction} algorithm does not add any colors to the multiset of colors used by $\mathcal{Q}$ as it builds $\mathcal{Q}_{\mbox{\em \tiny mountain}}$, then $\mathfrak{M}_{J}(\mathcal{Q}_{\mbox{\em \tiny mountain}}) \subseteq \mathfrak{M}_{J}(\mathcal{Q})$; since $\mathcal{Q}$ uses a minimum number of color $J$ edges, these two multisets of colors must be equal. 
The apex $\aelt$ of $\mathcal{Q}_{\mbox{\em \tiny mountain}}$ is an upper bound for both $\selt$ and $\telt$, so $\selt \vee \telt \leq \aelt$. 

Let $\mathcal{R}$ be the path that ascends from $\selt$ to $\selt \vee \telt$ along edges from $\mathcal{P}_{\mbox{\em \tiny mountain}}$, then ascends from $\selt \vee \telt $ up to $\aelt$, then descends from $\aelt$ to $\selt \vee \telt$ along the same edges, and then descends from $\selt \vee \telt$ to $\telt$ along edges from $\mathcal{P}_{\mbox{\em \tiny mountain}}$.
Observe that the `prefix' $\mathcal{Q}' := (\xelt_{0}=\selt,\ldots,\aelt)$ of the path $\mathcal{Q}_{\mbox{\em \tiny mountain}}$ is a path that ascends from $\selt$ to $\aelt$, and the `suffix' $\mathcal{Q}'' := (\aelt,\ldots,\xelt_{p}=\telt)$ descends from $\aelt$ to $\telt$. 
Similarly define the prefix $\mathcal{R}'$ and suffix $\mathcal{R}''$ of $\mathcal{R}$. 
By \TechnicalLemma, $\mathfrak{M}_{I}(\mathcal{Q}')=\mathfrak{M}_{I}(\mathcal{R}')$ and $\mathfrak{M}_{I}(\mathcal{Q}'')=\mathfrak{M}_{I}(\mathcal{R}'')$. 
So, $\mathfrak{M}_{I}(\mathcal{Q}_{\mbox{\em \tiny mountain}})=\mathfrak{M}_{I}(\mathcal{R})$. 
Of course, since $\mathcal{P}$ is obtained from $\mathcal{R}$ simply by omitting edges between $\selt \vee \telt$ and $\aelt$, then $\mathfrak{M}_{J}(\mathcal{P}) \subseteq \mathfrak{M}_{J}(\mathcal{Q}_{\mbox{\em \tiny mountain}})$. 
But, $\mathfrak{M}_{J}(\mathcal{Q}_{\mbox{\em \tiny mountain}}) = \mathfrak{M}_{J}(\mathcal{Q})$ is as small as possible. 
So, the multisets $\mathfrak{M}_{J}(\mathcal{P})$ and $\mathfrak{M}_{J}(\mathcal{Q})$ are the same.\hfill\QED 

Next we record some remarks relating to \MoveMinGameTheorem\ and offer comments to clarify some of our notation. 
In the next section, we apply \MoveMinGameTheorem\ to solve the Mixedmiddleswitch Puzzle.

\noindent 
{\bf \MoveMinSectionRemarks}

\vspace*{-0.2in}
\begin{enumerate}
\item[(1)]  The phrase `{\em God's algorithm}' is sometimes used in connection with move-minimizing puzzles that are modelled using Cayley graphs of related groups, perhaps most famously in connection with Rubik's Cube (see for example \cite{Joyner}). 
This phrase refers to an algorithm that minimizes the number of moves between two configurations of the puzzle, and `{\em God's number}' is the largest of all such move-minimizing numbers.  
In a similar sense, then, the preceding theorem determines God's number and God's algorithm for a move-minimizing puzzle between two specified objects whenever the puzzle digraph $\mathcal{G}$ is the covering digraph for a diamond-colored distributive lattice $L$. Specifically, part (1a) of \MoveMinGameTheorem\ says that God's number is the length of $L$, and part (1b) of \MoveMinGameTheorem\ prescribes God's algorithm. 
\item[(2)]  For modular lattices, we believe the results of the preceding theorem can be sharpened using Markowsky's theory of a poset of irreducibles for general lattices, see for example \cite{Markowsky}. 
\item[(3)]  At first glance, it may seem unlikely that any interesting puzzles could be modeled using diamond-colored modular/distributive lattices. 
In the next two sections, we present some move-minimizing puzzles that we believe serve as nice illustrations of this approach. 
We will explore, here and in subsequent papers, some of the many enumerative and algebraic contexts within which our move-minimizing puzzle digraphs seem naturally to occur. 
\item[(4)]  \MoveMinGameTheorem\ suggests a method for designing move-minimizing puzzles. 
Starting with a vertex-colored poset (or family of posets, if a puzzle is to be developed in a more generic way), one can try to devise an alternative coordinatization of the order ideals in such a way that the covering relations between order ideals in the associated diamond-colored distributive lattice have interesting, curious, and/or diversional descriptions in terms of these alternative coordinates. 
The move-minimizing puzzle one thusly obtains would use these alternative coordinates for order ideals as the set of objects and use some version of the (colored) covering relations as the set of  moves. 
\item[(5)] A note on notation: The names of the bijections presented in the remaining sections of the paper are given in a script-style font ({\tt {$\backslash$}mathscr$\{\cdot\}$} in {\tt LaTeX}), viz.\ $\mathscr{B}$ in the context of the Mixedmiddleswitch Puzzle and $\mathscr{M}$ for our recreational example in Section~\Closing. 
A typographical variation of `x' is used as the input variable for these functions (e.g.\ $\myspecialchar{x}$ for set-like coordinates, $x$ for a binary sequence, $\myx$ for a partition, $X$ for a columnar tableau, etc), and similar variations on `y' denote function output.  
Similar typographical variations of `s' and `t' are used to denote generic lattice elements in various coordinatizations. 
A calligraphic-style font ({\tt {$\backslash$}mathcal$\{\cdot\}$} in {\tt LaTeX}) together with the integer index $n$ is used to denote sets of objects in various  coordinates, viz.\ $\mathcal{Z}(n)$ and $\mathcal{B}(n)$ in the context of the Mixedmiddleswitch Puzzle, and $\mathcal{M}(n)$ and $\mathcal{C}(n)$ for our recreational puzzle in Section~\Closing. 
\end{enumerate}

\vspace{1ex} 
\noindent 
{\Large \bf \MixedMiddleSwitch\hspace*{0.15in}Solving the Mixedmiddleswitch Puzzle}

\noindent
To solve the Mixedmiddleswitch Puzzle on $n$ switches that was proposed in Section~\Intro, we require an explicit description of an isomorphism between the puzzle digraph and a diamond-colored modular or distributive lattice. 
Fix an integer $n \in (1,\infty)_{\mathbb{Z}}$. 
The objects comprising the vertex set of our puzzle digraph $\mathcal{B}(n)$ are the binary sequences of length $n$. 
We will codify the legal moves of the Mixedmiddleswitch Puzzle by creating the following colored and directed edges between binary sequences:   
For $i \in [1,n]_{\mathbb{Z}}$ we write $(s_1, s_2,\ldots, s_n) \myarrow{i} (t_1, t_2,\ldots, t_n)$ if and only if $s_{j}=t_{j}$ for each index $j \not= i$ and exactly one of the following holds: 
$1 < i < n$, $(s_{i-1}, s_{i}, s_{i+1}) = (0,0,1)$, and $(t_{i-1}, t_{i}, t_{i+1}) = (0,1,1)$; $1 < i < n$, $(s_{i-1}, s_{i}, s_{i+1}) = (1,1,0)$, and $(t_{i-1}, t_{i}, t_{i+1}) = (1,0,0)$; $i=1$, $(s_{1},s_{2}) = (0,1)$, and $(t_{1},t_{2}) = (1,1)$; $i=n$, $(s_{n-1},s_{n}) = (0,0)$, and $(t_{n-1},t_{n}) = (0,1)$; or $i=n$, $(s_{n-1},s_{n}) = (1,1)$, and $(t_{n-1},t_{n}) = (1,0)$. 
In such a case, we call $i$ the {\em color} of the corresponding mixedmiddleswitch move. 

Our strategy for solving the Mixedmiddleswitch Puzzle is to find a diamond-colored distributive lattice that is isomorphic to the puzzle digraph $\mathcal{B}(n)$.  The diamond-colored distributive lattice that effects this solution is the so-called minuscule lattice associated with the type $\myB_{n}$ root system / complex simple Lie algebra / Weyl group. 
We will have more to say about these algebraic contexts in the next section.   
For now, we simply describe the diamond-colored distributive lattice as a natural partial ordering of certain $n$-tuples that correspond to the subsets of $[1,n]_{\mathbb{Z}}$. 
For a lively and informative discussion of the role these distributive lattices played in solving an Erd\H{o}s problem, see R.\ A.\ Proctor's article \cite{PrMonthly}. 

Let $\mathcal{Z}(n)$ be the set 
\[\left\{\begin{array}{c}\mbox{integer sequences }\\ \myspecialchar{s} = (\myspecialchar{s}_{1}, \myspecialchar{s}_{2},\ldots, \myspecialchar{s}_{n})\end{array}\, \rule[-5mm]{0.25mm}{12mm}\, \begin{array}{c} n \geq \myspecialchar{s}_{1} \geq \cdots \geq \myspecialchar{s}_{n} \geq 0  \mbox{ where }\\ \myspecialchar{s}_{i} > \myspecialchar{s}_{i+1} \mbox{ if $i+1 \leq n$ and $\myspecialchar{s}_{i} \not= 0$}\end{array}\right\},\] 
partially ordered by component-wise comparison, i.e.\ $\myspecialchar{s} = (\myspecialchar{s}_{1}, \myspecialchar{s}_{2},\ldots, \myspecialchar{s}_{n}) \leq (\myspecialchar{t}_{1}, \myspecialchar{t}_{2},\ldots, \myspecialchar{t}_{n}) = \myspecialchar{t}$ for $\myspecialchar{s}, \myspecialchar{t} \in \mathcal{Z}(n)$ if and only if $\myspecialchar{s}_{i} \leq \myspecialchar{t}_{i}$ for all $i \in [1,n]_{\mathbb{Z}}$. 
One can see that $\myspecialchar{t}$ covers $\myspecialchar{s}$ in $\mathcal{Z}(n)$ if and only if there is some $k$ such that $\myspecialchar{t}_{k} = \myspecialchar{s}_{k}+1$ but $\myspecialchar{t}_{j} = \myspecialchar{s}_{j}$ for all indices $j \not= k$, in which case we give the covering digraph edge that corresponds to this covering relation the color $i := n+1-\myspecialchar{t}_{k}$ and write $\myspecialchar{s} \mylongarrow{i} \myspecialchar{t}$. 
It is easy to check that the resulting edge-colored poset is a diamond-colored distributive lattice. 
We refer to the $n$-tuples of $\mathcal{Z}(n)$ as {\em zero-cushioned $[1,n]$-subsets}. 
See \BFiveFigure\ for a depiction of $\mathcal{B}(n)$ and $\mathcal{Z}(n)$ when $n=5$. 

\begin{figure}[t]
\begin{center}
\setlength{\unitlength}{1.325cm}
\begin{picture}(6,10)
\put(-2,9.75){\large $\mathcal{Z}(5) \cong \mathcal{B}(5)$}
\put(0,0){\NEEdgeForBnLattice{Purple}{5}{0}{0}}
\put(1,0.6667){\NEEdgeForBnLattice{orange}{4}{0}{0}}
\put(2,1.3334){\NWEdgeForBnLattice{Green}{3}{0}{0}}
\put(2,1.3334){\NEEdgeForBnLattice{Purple}{5}{0}{0}}
\put(1,2){\NWEdgeForBnLattice{blue}{2}{0}{0}}
\put(1,2){\NEEdgeForBnLattice{Purple}{5}{0}{0}}
\put(3,2){\NWEdgeForBnLattice{Green}{3}{0}{0}}
\put(0,2.6667){\VerticalEdgeForBnLattice{red}{1}{0}{0}}
\put(0,2.6667){\NEEdgeForBnLattice{Purple}{5}{0}{0}}
\put(2,2.6667){\NWEdgeForBnLattice{blue}{2}{0}{0}}
\put(2,2.6667){\NEEdgeForBnLattice{orange}{4}{0}{0}}
\put(0,3.3334){\NEEdgeForBnLattice{Purple}{5}{0}{0}}
\put(1,3.3334){\VerticalEdgeForBnLattice{red}{1}{0}{0}}
\put(1,3.3334){\NEEdgeForBnLattice{orange}{4}{0}{0}}
\put(3,3.3334){\NWEdgeForBnLattice{blue}{2}{0}{0}}
\put(3,3.3334){\NEEdgeForBnLattice{Purple}{5}{0}{0}}
\put(1,4){\NEEdgeForBnLattice{orange}{4}{-0.25}{-0.15}}
\put(2,4){\NWEdgeForBnLattice{Green}{3}{-0.25}{0.15}}
\put(2,4){\VerticalEdgeForBnLattice{red}{1}{0}{0}}
\put(2,4){\NEEdgeForBnLattice{Purple}{5}{0}{0}}
\put(4,4){\NWEdgeForBnLattice{blue}{2}{0}{0}}
\put(1,4.6667){\VerticalEdgeForBnLattice{red}{1}{0}{0}}
\put(1,4.6667){\NEEdgeForBnLattice{Purple}{5}{0.25}{0.15}}
\put(2,4.6667){\NWEdgeForBnLattice{Green}{3}{0.25}{-0.15}}
\put(2,4.6667){\NEEdgeForBnLattice{Purple}{5}{-0.25}{-0.15}}
\put(3,4.6667){\NWEdgeForBnLattice{Green}{3}{-0.4}{0.15}}
\put(3,4.6667){\VerticalEdgeForBnLattice{red}{1}{0}{0}}
\put(1,5.3334){\NWEdgeForBnLattice{blue}{2}{0}{0}}
\put(1,5.3334){\NEEdgeForBnLattice{Purple}{5}{0}{0}}
\put(2,5.3334){\VerticalEdgeForBnLattice{red}{1}{0}{0}}
\put(2,5.3334){\NEEdgeForBnLattice{orange}{4}{0.25}{0.15}}
\put(3,5.3334){\NWEdgeForBnLattice{Green}{3}{0.25}{-0.15}}
\put(0,6){\NEEdgeForBnLattice{Purple}{5}{0}{0}}
\put(2,6){\NWEdgeForBnLattice{blue}{2}{0}{0}}
\put(2,6){\NEEdgeForBnLattice{orange}{4}{0}{0}}
\put(3,6){\VerticalEdgeForBnLattice{red}{1}{0}{0}}
\put(3,6){\NEEdgeForBnLattice{Purple}{5}{0}{0}}
\put(1,6.6667){\NEEdgeForBnLattice{orange}{4}{0}{0}}
\put(3,6.6667){\NWEdgeForBnLattice{blue}{2}{0}{0}}
\put(3,6.6667){\NEEdgeForBnLattice{Purple}{5}{0}{0}}
\put(4,6.6667){\VerticalEdgeForBnLattice{red}{1}{0}{0}}
\put(2,7.3334){\NWEdgeForBnLattice{Green}{3}{0}{0}}
\put(2,7.3334){\NEEdgeForBnLattice{Purple}{5}{0}{0}}
\put(4,7.3334){\NWEdgeForBnLattice{blue}{2}{0}{0}}
\put(1,8){\NEEdgeForBnLattice{Purple}{5}{0}{0}}
\put(3,8){\NWEdgeForBnLattice{Green}{3}{0}{0}}
\put(2,8.6667){\NEEdgeForBnLattice{orange}{4}{0}{0}}
\put(3,9.3334){\NEEdgeForBnLattice{Purple}{5}{0}{0}}
\put(0,0){\VertexForBnLattice{0}{0}{0}{0}{0}{0}{0}}
\put(0,0){\VertexForBnLatticeToo{0}{0}{0}{0}{0}{0}{0}}
\put(1,0.6667){\VertexForBnLattice{1}{0}{0}{0}{0}{0}{0}}
\put(1,0.6667){\VertexForBnLatticeToo{0}{0}{0}{0}{1}{0}{0}}
\put(2,1.3334){\VertexForBnLattice{2}{0}{0}{0}{0}{0}{0}}
\put(2,1.3334){\VertexForBnLatticeToo{0}{0}{0}{1}{1}{0}{0}}
\put(1,2){\VertexForBnLattice{3}{0}{0}{0}{0}{-1.05}{0}}
\put(1,2){\VertexForBnLatticeToo{0}{0}{1}{1}{1}{-1.05}{0}}
\put(3,2){\VertexForBnLattice{2}{1}{0}{0}{0}{-0.1}{0}}
\put(3,2){\VertexForBnLatticeToo{0}{0}{0}{1}{0}{-0.1}{0}}
\put(0,2.6667){\VertexForBnLattice{4}{0}{0}{0}{0}{-0.9}{0}}
\put(0,2.6667){\VertexForBnLatticeToo{0}{1}{1}{1}{1}{-0.9}{0}}
\put(2,2.6667){\VertexForBnLattice{3}{1}{0}{0}{0}{0}{0}}
\put(2,2.6667){\VertexForBnLatticeToo{0}{0}{1}{1}{0}{0}{0}}
\put(0,3.3334){\VertexForBnLattice{5}{0}{0}{0}{0}{-0.9}{0}}
\put(0,3.3334){\VertexForBnLatticeToo{1}{1}{1}{1}{1}{-0.9}{0}}
\put(1,3.3334){\VertexForBnLattice{4}{1}{0}{0}{0}{0}{0}}
\put(1,3.3334){\VertexForBnLatticeToo{0}{1}{1}{1}{0}{0}{0}}
\put(3,3.3334){\VertexForBnLattice{3}{2}{0}{0}{0}{0}{0}}
\put(3,3.3334){\VertexForBnLatticeToo{0}{0}{1}{0}{0}{0}{0}}
\put(1,4){\VertexForBnLattice{5}{1}{0}{0}{0}{-1}{0}}
\put(1,4){\VertexForBnLatticeToo{1}{1}{1}{1}{0}{-1}{0}}
\put(2,4){\VertexForBnLattice{4}{2}{0}{0}{0}{0}{0}}
\put(2,4){\VertexForBnLatticeToo{0}{1}{1}{0}{0}{0}{0}}
\put(4,4){\VertexForBnLattice{3}{2}{1}{0}{0}{-0.1}{0}}
\put(4,4){\VertexForBnLatticeToo{0}{0}{1}{0}{1}{-0.1}{0}}
\put(1,4.6667){\VertexForBnLattice{4}{3}{0}{0}{0}{-1}{0}}
\put(1,4.6667){\VertexForBnLatticeToo{0}{1}{0}{0}{0}{-1}{0}}
\put(2,4.6667){\VertexForBnLattice{5}{2}{0}{0}{0}{0}{0}}
\put(2,4.6667){\VertexForBnLatticeToo{1}{1}{1}{0}{0}{0}{0}}
\put(3,4.6667){\VertexForBnLattice{4}{2}{1}{0}{0}{0}{0}}
\put(3,4.6667){\VertexForBnLatticeToo{0}{1}{1}{0}{1}{0}{0}}
\put(1,5.3334){\VertexForBnLattice{5}{3}{0}{0}{0}{-1}{0}}
\put(1,5.3334){\VertexForBnLatticeToo{1}{1}{0}{0}{0}{-1}{0}}
\put(2,5.3334){\VertexForBnLattice{4}{3}{1}{0}{0}{0}{0}}
\put(2,5.3334){\VertexForBnLatticeToo{0}{1}{0}{0}{1}{0}{0}}
\put(3,5.3334){\VertexForBnLattice{5}{2}{1}{0}{0}{0}{0}}
\put(3,5.3334){\VertexForBnLatticeToo{1}{1}{1}{0}{1}{0}{0}}
\put(0,6){\VertexForBnLattice{5}{4}{0}{0}{0}{-1}{0}}
\put(0,6){\VertexForBnLatticeToo{1}{0}{0}{0}{0}{-1}{0}}
\put(2,6){\VertexForBnLattice{5}{3}{1}{0}{0}{0}{0}}
\put(2,6){\VertexForBnLatticeToo{1}{1}{0}{0}{1}{0}{0}}
\put(3,6){\VertexForBnLattice{4}{3}{2}{0}{0}{0}{0}}
\put(3,6){\VertexForBnLatticeToo{0}{1}{0}{1}{1}{0}{0}}
\put(1,6.6667){\VertexForBnLattice{5}{4}{1}{0}{0}{-1}{0}}
\put(1,6.6667){\VertexForBnLatticeToo{1}{0}{0}{0}{1}{-1}{0}}
\put(3,6.6667){\VertexForBnLattice{5}{3}{2}{0}{0}{-1}{0}}
\put(3,6.6667){\VertexForBnLatticeToo{1}{1}{0}{1}{1}{-1}{0}}
\put(4,6.6667){\VertexForBnLattice{4}{3}{2}{1}{0}{0}{0}}
\put(4,6.6667){\VertexForBnLatticeToo{0}{1}{0}{1}{0}{0}{0}}
\put(2,7.3334){\VertexForBnLattice{5}{4}{2}{0}{0}{-1}{0}}
\put(2,7.3334){\VertexForBnLatticeToo{1}{0}{0}{1}{1}{-1}{0}}
\put(4,7.3334){\VertexForBnLattice{5}{3}{2}{1}{0}{0}{0}}
\put(4,7.3334){\VertexForBnLatticeToo{1}{1}{0}{1}{0}{0}{0}}
\put(1,8){\VertexForBnLattice{5}{4}{3}{0}{0}{-1}{0}}
\put(1,8){\VertexForBnLatticeToo{1}{0}{1}{1}{1}{-1}{0}}
\put(3,8){\VertexForBnLattice{5}{4}{2}{1}{0}{0}{0}}
\put(3,8){\VertexForBnLatticeToo{1}{0}{0}{1}{0}{0}{0}}
\put(2,8.6667){\VertexForBnLattice{5}{4}{3}{1}{0}{-1}{0}}
\put(2,8.6667){\VertexForBnLatticeToo{1}{0}{1}{1}{0}{-1}{0}}
\put(3,9.3334){\VertexForBnLattice{5}{4}{3}{2}{0}{-1}{0}}
\put(3,9.3334){\VertexForBnLatticeToo{1}{0}{1}{0}{0}{-1}{0}}
\put(4,10){\VertexForBnLattice{5}{4}{3}{2}{1}{-1}{0}}
\put(4,10){\VertexForBnLatticeToo{1}{0}{1}{0}{1}{-1}{0}}
\put(-3.1,-0.5){\small {\bf \BFiveFigure}\ \  The\hspace*{0.2mm} diamond-colored\hspace*{0.2mm} distributive\hspace*{0.2mm} lattice\hspace*{0.2mm} $\mathcal{Z}(n) \cong \mathcal{B}(n)$\hspace*{0.2mm} with\hspace*{0.2mm} $n=5$.\ \  At\hspace*{0.2mm} each\hspace*{0.2mm} vertex,\hspace*{0.2mm} the\hspace*{0.2mm} top} 
\put(-3.1,-0.8){\small coordinates are for the associated zero-cushioned $[1,5]$-subset; underneath the zero-cushioned $[1,5]$-sub-}
\put(-3.1,-1.1){\small set is the associated mixedmiddleswitch binary sequence.}
\end{picture}
\end{center}

\vspace*{0.75cm}
\hspace*{1in}$\overline{\hspace*{4.5in}}$
\end{figure}

We now describe an isomorphism $\mathscr{B}: \mathcal{Z}(n) \longrightarrow \mathcal{B}(n)$. 
Let $\myspecialchar{x} = (\myspecialchar{x}_{1}, \myspecialchar{x}_{2},\ldots, \myspecialchar{x}_{n})$ be an element of $\mathcal{Z}(n)$. 
For purposes of the next computation, let $\myspecialchar{x}_{0} := n$ and $\myspecialchar{x}_{n+1} := 0$.  
Informally, the binary sequence $y=\mathscr{B}(\myspecialchar{x})$\ \  will alternate between strings of $0$'s and strings of $1$'s in order to represent nonconsecutive entries of $\myspecialchar{x}$. 
More precisely, define a sequence of sets $\mathcal{I}(\myspecialchar{x}) := (\mathcal{I}_{0}, \mathcal{I}_{1}, \ldots, \mathcal{I}_{n})$ by the rule $\mathcal{I}_{i} := [n+1-\myspecialchar{x}_{i}, n-\myspecialchar{x}_{i+1}]_{\mathbb{Z}}$ for each $i \in [0,n]_{\mathbb{Z}}$. 
So $\mathcal{I}_{i}$ is a sequence of consecutive integers (to be thought of as indices) and is empty if $n+1-\myspecialchar{x}_{i} > n-\myspecialchar{x}_{i+1}$. 
Note that $\mathcal{I}_{i_1} \cap \mathcal{I}_{i_2} = \emptyset$ if $i_{1} \not= i_{2}$ and that $\displaystyle \mbox{\Large $\cup$}_{i=0}^{n}\, \mathcal{I}_{i} = \{1, 2, \ldots, n\}$. 
For each $j \in [1,n]_{\mathbb{Z}}$, set $y_{j} := i\, (\mbox{mod $2$})$ if $j \in \mathcal{I}_{i}$. 
Let $y = (y_{1}, y_{2}, \ldots, y_{n})$, and declare that $\mathscr{B}(\myspecialchar{x}) := y$.
We reverse this process by defining a function $\mathscr{B}': \mathcal{B}(n) \longrightarrow \mathcal{Z}(n)$ as follows. 
Let $y = (y_{1}, \ldots, y_{n}) \in \mathcal{B}(n)$, and set $y_{0} := 0$ and $y_{n+1} := 0$. 
Define the set \[\change(y) := \bigg\{j \in [1,n]_{\mathbb{Z}}\, \rule[-3mm]{0.2mm}{7.5mm}\ y_{j-1} \not= y_{j}\bigg\},\] which is the set of indices where the binary sequence $y$ changes from $0$ to $1$ (or vice-versa) as we read from left to right. 
Let $k := |\change(y)|$, and write the elements $j_{1}, \ldots j_{k}$ of $\change(y)$ in sequence so that $1 \leq j_{1} < j_{2} < \cdots < j_{k} \leq n$. 
(Of course, this sequence is empty if $k=0$.) 
For each $i \in [1,n]_{\mathbb{Z}}$, define 
\[\myspecialchar{x}_{i} := \left\{\begin{array}{cl}
n+1-j_{i} & \hspace*{0.1in}\mbox{if } 1 \leq i \leq k\\
0 & \hspace*{0.1in}\mbox{if } k+1 \leq i \leq n
\end{array}\right.\]
and let $\myspecialchar{x} := (\myspecialchar{x}_{1}, \myspecialchar{x}_{2},\ldots, \myspecialchar{x}_{n})$. 
Set $\mathscr{B}'(y) := \myspecialchar{x}$. 
The following additional construction will be useful in understanding $\mathscr{B}'$ as the inverse of $\mathscr{B}$. 
We set $j_{0} := 1$ and $j_{k+1} := n+1$. 
Now define a sequence of sets $\mathcal{J}(y) := (\mathcal{J}_{0}, \mathcal{J}_{1}, \ldots, \mathcal{J}_{n})$ by the rule 
\[\mathcal{J}_{i} := \left\{\begin{array}{cl}
[j_{i}, j_{i+1}-1]_{\mathbb{Z}} & \hspace*{0.1in}\mbox{if } i \in [1,k]_{\mathbb{Z}}\\
\emptyset & \hspace*{0.1in}\mbox{if } i \in [k+1,n]_{\mathbb{Z}}
\end{array}\right.\]
Note that $\mathcal{J}_{i_1} \cap \mathcal{J}_{i_2} = \emptyset$ if $i_{1} \not= i_{2}$ and that $\displaystyle \mbox{\Large $\cup$}_{i=0}^{n}\, \mathcal{J}_{i} = [1,n]_{\mathbb{Z}}$. 
Observe now that $\mathcal{I}(\myspecialchar{x}) = \mathcal{J}(\mathscr{B}(\myspecialchar{x}))$ and that $\mathcal{J}(y) = \mathcal{I}(\mathscr{B}'(y))$. 
It follows that $\mathscr{B}' = \mathscr{B}^{-1}$. 

\noindent
{\bf \TypeBLatticeIsoProp}\ \ 
{\sl The above set mapping $\mathscr{B}$ is an isomorphism of edge-colored directed graphs, so we have $\mathcal{Z}(n) \cong \mathcal{B}(n)$. In particular, each of these diamond-colored distributive lattices has $2^{n}$ elements and length} {\mbox{\footnotesize $\bigg(\!\!\!\begin{array}{c}\mbox{\footnotesize $n$}\\ \mbox{\footnotesize $2$}\end{array}\!\!\!\bigg)$}}.

{\em Proof.} The latter claims about the order and length of $\mathcal{Z}(n)$ are routine, so we only need to prove that $\mathscr{B}$ is an isomorphism. 
We have already seen that $\mathscr{B}$ is a bijection whose inverse is $\mathscr{B}^{-1} = \mathscr{B}'$. 
It only remains to be checked that $\mathscr{B}$ and $\mathscr{B}^{-1}$ are edge- and edge-color-preserving. 

So suppose $\myspecialchar{s} = (\myspecialchar{s}_{1}, \myspecialchar{s}_{2},\ldots, \myspecialchar{s}_{n})$ and $\myspecialchar{t} = (\myspecialchar{t}_{1}, \myspecialchar{t}_{2},\ldots, \myspecialchar{t}_{n})$ are in $\mathcal{Z}(n)$ with $\myspecialchar{s} \myarrow{i} \myspecialchar{t}$.  
Then for some $p \in [1,n]_{\mathbb{Z}}$, we have $\myspecialchar{t}_{p} = \myspecialchar{s}_{p}+1$  while $\myspecialchar{t}_{q}  = \myspecialchar{s}_{q}$ when $q \not= p$, in which case $i = n+1-\myspecialchar{t}_{p} = n-\myspecialchar{s}_{p}$.  

First, consider the case that $i=1$. 
Note that $p=1$ with  $\myspecialchar{t}_{1} = n$ and $\myspecialchar{s}_{1} = n-1$.  
Moreover $\myspecialchar{t}_{2} = \myspecialchar{s}_{2} < n-1 = \myspecialchar{s}_{1}$. 
In particular, the latter inequality implies that $n-\myspecialchar{t}_{2} = n-\myspecialchar{s}_{2} \geq 2$. 
Thus, we have
\[\mathcal{I}_{0}(\myspecialchar{t}) = \emptyset \mbox{\ \ \ and\ \ \ } \mathcal{I}_{1}(\myspecialchar{t}) = [1, n-\myspecialchar{t}_{2}]_{\mathbb{Z}}\]
while 
\[\mathcal{I}_{0}(\myspecialchar{s}) = [1, 1]_{\mathbb{Z}} \mbox{\ \ \ and\ \ \ } \mathcal{I}_{1}(\myspecialchar{s}) = [2, n-\myspecialchar{s}_{2}]_{\mathbb{Z}}.\]
Then the $1^{\mbox{\footnotesize st}}$ and $2^{\mbox{\footnotesize nd}}$ coordinates of $\mathscr{B}(\myspecialchar{t})$ are $(1,1)$, and coordinates in the same positions of $\mathscr{B}(\myspecialchar{s})$ are $(0,1)$. 
Otherwise, the coordinates of $\mathscr{B}(\myspecialchar{s})$ and $\mathscr{B}(\myspecialchar{t})$ agree. 
Therefore $\mathscr{B}(\myspecialchar{s}) \myarrow{1} \mathscr{B}(\myspecialchar{t})$ in $\mathcal{B}(n)$. 

Second, consider the case that $i=n$. 
Since $i = n+1-\myspecialchar{t}_{p} = n-\myspecialchar{s}_{p}$, we have $\myspecialchar{t}_{p} = 1$ and $\myspecialchar{s}_{p} = 0$. 
Moreover, $\myspecialchar{s}_{p-1} = \myspecialchar{t}_{p-1} \geq 2$, so $n+1-\myspecialchar{s}_{p-1} = n+1-\myspecialchar{t}_{p-1} \leq n-1$. 
Observe that  
\[\mathcal{I}_{p-1}(\myspecialchar{t}) = [n+1-\myspecialchar{t}_{p-1},n-1]_{\mathbb{Z}} \mbox{\ \ \ and\ \ \ } \mathcal{I}_{p}(\myspecialchar{t}) = [n,n]_{\mathbb{Z}}\]
while 
\[\mathcal{I}_{p-1}(\myspecialchar{s}) = [n+1-\myspecialchar{s}_{p-1},n]_{\mathbb{Z}} \mbox{\ \ \ and\ \ \ } \mathcal{I}_{p}(\myspecialchar{s}) = \emptyset.\] 
Then the $(n-1)^{\mbox{\footnotesize st}}$ and $n^{\mbox{\footnotesize th}}$ coordinates of $\mathscr{B}(\myspecialchar{t})$ are $(0,1)$ (if $p$ is odd) and $(1,0)$ (if $p$ is even), and coordinates in the same positions of $\mathscr{B}(\myspecialchar{s})$ are $(0,0)$ when $p$ is odd  and $(1,1)$ when $p$ is even. 
Otherwise, the coordinates of $\mathscr{B}(\myspecialchar{s})$ and $\mathscr{B}(\myspecialchar{t})$ agree. 
Therefore $\mathscr{B}(\myspecialchar{s}) \myarrow{n} \mathscr{B}(\myspecialchar{t})$ in $\mathcal{B}(n)$. 

Third (and last), consider the case that $1 < i < n$. 
The $(p-1)^{\mbox{\footnotesize st}}$ and $p^{\mbox{\footnotesize th}}$ integer intervals of  $\mathcal{I}(\myspecialchar{t})$ and $\mathcal{I}(\myspecialchar{s})$ are respectively 
\[\mathcal{I}_{p-1}(\myspecialchar{t}) = [n+1-\myspecialchar{t}_{p-1}, n-\myspecialchar{t}_{p}]_{\mathbb{Z}} \mbox{\ \ \ and\ \ \ } \mathcal{I}_{p}(\myspecialchar{t}) = [n+1-\myspecialchar{t}_{p}, n-\myspecialchar{t}_{p+1}]_{\mathbb{Z}}\]
with  
\[\mathcal{I}_{p-1}(\myspecialchar{s}) = [n+1-\myspecialchar{s}_{p-1}, n-\myspecialchar{s}_{p}]_{\mathbb{Z}} \mbox{\ \ \ and\ \ \ } \mathcal{I}_{p}(\myspecialchar{s}) = [n+1-\myspecialchar{s}_{p}, n-\myspecialchar{s}_{p+1}]_{\mathbb{Z}}.\]

Since $i = n+1-\myspecialchar{t}_{p} = n-\myspecialchar{s}_{p}$, then $i-1= n-\myspecialchar{t}_{p} = n-1-\myspecialchar{s}_{p}$ and $i+1 = n+2-\myspecialchar{t}_{p} = n+1-\myspecialchar{s}_{p}$. 
Also, $n+1-\myspecialchar{t}_{p-1} = n+1-\myspecialchar{s}_{p-1} \leq i-1$ and $n-\myspecialchar{s}_{p+1} = n-\myspecialchar{t}_{p+1} \geq i+1$. 
Let us think of the preceding integer intervals as 
\[\mathcal{I}_{p-1}(\myspecialchar{t}) = [\mbox{``$\leq i-1$''}, i-1]_{\mathbb{Z}} \mbox{\ \ \ and\ \ \ } \mathcal{I}_{p}(\myspecialchar{t}) = [i, \mbox{``$\geq i+1$''}]_{\mathbb{Z}}\] 
with 
\[\mathcal{I}_{p-1}(\myspecialchar{s}) = [\mbox{``$\leq i-1$''}, i]_{\mathbb{Z}} \mbox{\ \ \ and\ \ \ } \mathcal{I}_{p}(\myspecialchar{s}) = [i+1, \mbox{``$\geq i+1$''}]_{\mathbb{Z}}.\]
When $p$ is even, the $(i-1)^{\mbox{\footnotesize st}}$, $i^{\mbox{\footnotesize th}}$, and $(i+1)^{\mbox{\footnotesize st}}$ coordinates of $\mathscr{B}(\myspecialchar{t})$ are $(1,0,0)$ and coordinates in the same positions of $\mathscr{B}(\myspecialchar{s})$ are $(1,1,0)$;  
when $p$ is odd, coordinates in these positions are $(0,1,1)$ for $\mathscr{B}(\myspecialchar{t})$ and$(0,0,1)$ for $\mathscr{B}(\myspecialchar{s})$. 
Otherwise, the coordinates of $\mathscr{B}(\myspecialchar{s})$ and $\mathscr{B}(\myspecialchar{t})$ agree. 
Therefore $\mathscr{B}(\myspecialchar{s}) \myarrow{i} \mathscr{B}(\myspecialchar{t})$ in $\mathcal{B}(n)$. 

Thus $\mathscr{B}$ preserves edges and edge colors. Similar reasoning can be used see that $\mathscr{B}^{-1}$ also preserves edges and edge colors.\hfill\QED

In view of \TypeBLatticeIsoProp, \MoveMinGameTheorem\ applies to the Mixedmiddleswitch Puzzle. 
A complete solution to the puzzle would utilize $\mathbf{j}_{\mbox{\em \scriptsize color}}(\mathcal{Z}(n))$ to express shortest paths between objects of $\mathcal{B}(n)$ as sequences of elements from $\mathbf{j}_{\mbox{\em \scriptsize color}}(\mathcal{Z}(n))$. 
This latter vertex-colored poset is well-known as the type $\myB_{n}$ minuscule poset, cf.\ \cite{PrEur}. 
The next result offers a partial solution to the Mixedmiddleswitch Puzzle by explicitly addressing move-minimizing puzzle questions (1a) and (2a). 
For the purpose of stating this result, we make the following definitions: Let $\myspecialchar{x} = (\myspecialchar{x}_{1}, \myspecialchar{x}_{2},\ldots, \myspecialchar{x}_{n})$ and $\myspecialchar{y} = (\myspecialchar{y}_{1}, \myspecialchar{y}_{2},\ldots, \myspecialchar{y}_{n})$ be any zero-cushioned $[1,n]$-subsets with $\myspecialchar{x} \leq \myspecialchar{y}$. 
Let $C(\myspecialchar{x},\myspecialchar{y})$ denote the set $\left\{j \in [1,n]_{\mathbb{Z}}\  \rule[-1.5mm]{0.2mm}{5mm}\  n+1-\myspecialchar{x}_{i} > j \geq n+1-\myspecialchar{y}_{i} \mbox{ for some } i \in [1,n]_{\mathbb{Z}}\right\}$; and, for any $j \in [1,n]_{\mathbb{Z}}$, let $\#_{j}(\myspecialchar{x},\myspecialchar{y}) := \rule[-1.5mm]{0.2mm}{5mm}\{i \in [1,n]_{\mathbb{Z}}\ |\, n+1-\myspecialchar{x}_{i} > j \geq n+1-\myspecialchar{y}_{i}\}\rule[-1.5mm]{0.2mm}{5mm}$. 

\noindent
{\bf \MixedMiddleSwitchProp}\ \ 
{\sl Let $s$ and $t$ be binary sequences in $\mathcal{B}(n)$, thought of as the starting and terminating objects for Mixedmiddleswitch Puzzle play, and let $\myspecialchar{s} = (\myspecialchar{s}_{1}, \myspecialchar{s}_{2},\ldots, \myspecialchar{s}_{n}) := \mathscr{B}^{-1}(s)$ and $\myspecialchar{t} = (\myspecialchar{t}_{1}, \myspecialchar{t}_{2},\ldots, \myspecialchar{t}_{n}) := \mathscr{B}^{-1}(t)$ be their preimages in $\mathcal{Z}(n)$ under our edge-colored digraph isomorphism $\mathscr{B}$.  
Let $\rho_{_\mathcal{Z}}$ denote the rank function for $\mathcal{Z}(n)$. Then:} 
\begin{enumerate}
\item[{\sl (1)}] {\sl We have $\rho_{_\mathcal{Z}}(\myspecialchar{s})=\sum \myspecialchar{s}_i$ and $\rho_{_\mathcal{Z}}(\myspecialchar{t})=\sum \myspecialchar{t}_i$.  
Also, $\myspecialchar{s} \vee \myspecialchar{t} = (\mathsf{max}(\myspecialchar{s}_i,\myspecialchar{t}_i))_{i=1,\ldots,n}$ and $\myspecialchar{s} \wedge \myspecialchar{t} = (\mathsf{min}(\myspecialchar{s}_i,\myspecialchar{t}_i))_{i=1,\ldots,n}$. 
That is, the join and meet of  $\myspecialchar{s}$ and $\myspecialchar{t}$ are just the component-wise max and component-wise min, respectively.}  
{\sl Moreover, we have} 
\[\mathsf{dist}(s,t)\ =\ \sum_{i=1}^{n}\left(\rule[-1.5mm]{0mm}{5mm}2\mathsf{max}(\myspecialchar{s}_i,\myspecialchar{t}_i) - \myspecialchar{s}_i - \myspecialchar{t}_i\rule[-1.5mm]{0mm}{5mm}\right)\ =\ \sum_{i=1}^{n}\left(\rule[-1.5mm]{0mm}{5mm}\myspecialchar{s}_i + \myspecialchar{t}_i - 2\mathsf{min}(\myspecialchar{s}_i,\myspecialchar{t}_i)\rule[-1.5mm]{0mm}{5mm}\right)\ \leq\ \mbox{\footnotesize $\bigg(\!\!\!\begin{array}{c}\mbox{\footnotesize $n$}\\ \mbox{\footnotesize $2$}\end{array}\!\!\!\bigg)$}\ <\ \infty.\]
{\sl In particular, the preceding calculation of $\mathsf{dist}(s,t)$ explicitly solves question (1a) for the Mixedmiddleswitch Puzzle.}
\item[{\sl (2)}] {\sl Keep the notation of part (1), and let $J \subseteq \{1,2,\ldots,n\}$.  
Question (2a) for the Mixedmiddleswitch Puzzle is solved as follows. 
The minimum number of mixedmiddleswitch moves with colors from the set $J$ that are required in order to obtain $t$ from $s$ is the left-hand-side quantity in the following inequality:} 
\[\sum_{j \in J}\left(\#_{j}(\myspecialchar{s},\myspecialchar{s}\vee\myspecialchar{t})\rule[-1.5mm]{0mm}{5mm}+\rule[-1.5mm]{0mm}{5mm}\#_{j}(\myspecialchar{t},\myspecialchar{s}\vee\myspecialchar{t})\right)\ \leq\ \mathsf{dist}(s,t).\] 
{\sl This minimum count of color $J$ mixedmiddlswitch moves is attained by any shortest sequence of mixedmiddleswitch moves used to maneuver from $s$ to $t$. 
Moreover, $t$ can be obtained from $s$ via mixedmiddleswitch moves with colors only from the set $J$ if and only if $C(\myspecialchar{s},\myspecialchar{s}\vee\myspecialchar{t}) \cup C(\myspecialchar{t},\myspecialchar{s}\vee\myspecialchar{t}) \subseteq J$ if and only if} 
\[\mathsf{dist}(s,t)\ =\ \sum_{j \in J}\left(\#_{j}(\myspecialchar{s},\myspecialchar{s}\vee\myspecialchar{t})\rule[-1.5mm]{0mm}{5mm}+\rule[-1.5mm]{0mm}{5mm}\#_{j}(\myspecialchar{t},\myspecialchar{s}\vee\myspecialchar{t})\right).\]
{\sl All the statements of part (2) remain valid if we appropriately replace each occurrence of `$\myspecialchar{s}\vee\myspecialchar{t}$' with `$\myspecialchar{s}\wedge\myspecialchar{t}$'.}
\end{enumerate}

{\em Proof.} \MoveMinGameTheorem\ applies here, as \TypeBLatticeIsoProp\ has established that the puzzle digraph $\mathcal{B}(n)$ is isomorphic to the distributive lattice $\mathcal{Z}(n)$. 
To complete the proof of \MixedMiddleSwitchProp, we must simply verify that the meet, join, lattice length, element rank, and path length computations in $\mathcal{Z}(n)$ are correctly applied here. 
The clams in part {\sl (1)} of the proposition statement concerning meet, join, and rank computations in $\mathcal{Z}(n)$ follow from the definitions; that $\mathcal{Z}(n)$ has length $\mbox{\footnotesize $\bigg(\!\!\!\begin{array}{c}\mbox{\footnotesize $n$}\\ \mbox{\footnotesize $2$}\end{array}\!\!\!\bigg)$}$ was observed in \TypeBLatticeIsoProp. 
So, the displayed formulas and inequality in part {\sl (1)} of the proposition statement follow from \MoveMinGameTheorem.1. 

Now consider a subset $J$ of $\{1,2,\ldots,n\}$. 
Recall that $\myspecialchar{p} \myarrow{j} \myspecialchar{q}$ in $\mathcal{Z}(n)$ for some $j \in J$ if and only if for some $i \in [1,n]_{\mathbb{Z}}$ we have $n+1-\myspecialchar{p}_{i} - 1 = j = n+1-\myspecialchar{q}_{i}$ while $\myspecialchar{p}_{k} = \myspecialchar{q}_{k}$ for each $k \in [1,n]_{\mathbb{Z}} \setminus \{i\}$. 
From this observation, it follows that if $\myspecialchar{x} \leq \myspecialchar{y}$ in $\mathcal{Z}(n)$, then $C(\myspecialchar{x},\myspecialchar{y}) \cap J$ is exactly the set of color $J$ edges in any path in $\mathcal{Z}(n)$ from $\myspecialchar{x}$ up to $\myspecialchar{y}$ and $\sum_{j \in J}\left(\#_{j}(\myspecialchar{x},\myspecialchar{y})\rule[-1.5mm]{0mm}{5mm}\right)$ is the total number of color $J$ edges encountered in such a path. 
The claims of part {\sl (2)} of the proposition statement now follow from \MoveMinGameTheorem.2.\hfill\QED

So, for the question posed in Section~\Intro, the minimum number of mixedmiddleswitch moves needed to maneuver from $s := (0,0,0,0,0)$ to $t := (0,1,0,1,0)$ is the same as the distance from $\myspecialchar{s} = (0,0,0,0,0)$ to $\myspecialchar{t} = (4,3,2,1,0)$ in $\mathcal{Z}(5)$. 
Since $\myspecialchar{s} = \mathsf{min}(\mathcal{Z}(5))$ and $\rho_{_\mathcal{Z}}(\myspecialchar{t}) = 10$, then $\dist(\myspecialchar{s},\myspecialchar{t}) = 10$, so the minimum number of mixedmiddleswitch moves from $(0,0,0,0,0)$ to $(0,1,0,1,0)$ is 10.

\vspace{1ex} 
\noindent 
{\Large \bf \Closing\hspace*{0.15in}A recreational move-minimizing puzzle}

\noindent 
We close with a challenge, which is to solve a move-minimizing puzzle we call `Snakes on a Ming'antu Square' developed by the first- and third-listed authors. 
(We give an ample hint towards the end of this section.)  
Here are some of the basics of our set-up. 
Let $n$ be a positive integer, and fix once and for all an $n \times n$ board of squares. 
The squares of the board are coordinatized by integer pairs $(i,j)$ in the way we normally index the entries of square matrices; that is, square $(i,j)$ is in row $i$ and column $j$ of our board, where $i,j \in [1,n]_{\mathbb{Z}}$, and the main diagonal consists of the pairs $(i,i)$. 
Movement `left' or `west' is in the direction of decreasing $j$ along a given row $i$, movement `down' or `south' is in the direction of increasing $i$ within a given column $j$, etc.   Puzzle play will require us sometimes to cover and sometimes to uncover a given square with a $1 \times 1$ `tile'.

The objects of our puzzle are what we call `Ming'antu tilings' of our board and comprise a set $\mathcal{M}(n)$. 
A {\em Ming'antu tiling} is a placement of tiles on our $n \times n$ game board satisfying the following rules: 
(1) If a square $(i,j)$ is covered by a tile, then so are all squares above or to the left of $(i,j)$; 
and (2) If a main diagonal square $(i,i)$ is covered with a tile, then the number of tiled squares below $(i,i)$ in column $i$ does not exceed the number of tiled squares to the  right of $(i,i)$ in row $i$. 

The puzzle allows a player to add or remove a set of tiles from a given Ming'antu tiling if the underlying squares form a certain kind of snake-like pattern.  Specifically, a {\em centered southwesterly snake} is a sequence of $m$ squares chosen from our board (where $m$ is some positive integer) such that (1) the sequence begins at some fixed square and from any square in the sequence the succeeding square is either one step to its South or to its West and (2) the $\lfloor \frac{m+1}{2} \rfloor^{\mbox{\footnotesize st}}$ square in the sequence is on the main diagonal of the game board.  Given a Ming'antu tiling $\myx$, a {\em centered southwesterly snake move} is either the addition of tiles to some centered southwesterly snake whose squares are not tiled by $\myx$  or else the removal of tiles from $\myx$ along some centered southwesterly snake whose squares are amongst those tiled by $\myx$; the move is legal only if the result is another Ming'antu tiling. 

\renewcommand{\thefootnote}{9}
{\sl Snakes on a Ming'antu Square} is the move-minimizing puzzle in which the player is given two Ming'antu tilings of a square board and must move from one tiling to the other by (legal) centered southwesterly snake moves (hereafter called `snake moves'). 
Our challenge to the reader is to solve the associated move-minimizing puzzle.  
Here's a specific instance of the puzzle that might serve as a helpful warm-up: 
What is the minimum number of snake moves required to move from $\mys$ to $\myt$, where these Ming'antu tilings are identified as in \ChallengeFigure?\footnote{The reader may contact one of the authors for an answer to this question.}

\begin{figure}[h]
\begin{center}
\setlength{\unitlength}{0.4cm}
\begin{picture}(6,4)
\thicklines
\put(-6.5,1.6){\Large $\mylargers = $}
\put(-4.5,0){\multiput(0,3)(1,0){4}{\multiput(0.05,0)(0.025,0){39}{\textcolor{gray}{\line(0,1){1}}}}
\multiput(0,2)(1,0){4}{\multiput(0.05,0)(0.025,0){39}{\textcolor{gray}{\line(0,1){1}}}}
\multiput(0,1)(1,0){1}{\multiput(0.05,0)(0.025,0){39}{\textcolor{gray}{\line(0,1){1}}}}
\multiput(0,0)(1,0){5}{\line(0,1){4}}
\multiput(0,0)(0,1){5}{\line(1,0){4}}}
\put(0.8,1.6){\huge $\stackrel{\mbox{\tiny how?}}{\rightsquigarrow}$}
\put(8,1.6){\Large$ = \mylargert$}
\put(3.5,0){\multiput(0,3)(1,0){1}{\multiput(0.05,0)(0.025,0){39}{\textcolor{gray}{\line(0,1){1}}}}
\multiput(0,0)(1,0){5}{\line(0,1){4}}
\multiput(0,0)(0,1){5}{\line(1,0){4}}}
\put(-8.6,-1.75){\small {\bf \ChallengeFigure}\ \ A Snakes-on-a-Ming'antu-Square challenge!}
\end{picture}
\end{center}

\vspace*{0.0cm}
\hspace*{4.25cm}$\overline{\hspace*{7cm}}$
\end{figure}

Here is the promised hint. 
We will codify snake moves using colored and directed edges as follows:
Write $\mys \myarrow{i} \myt$ for Ming'antu tilings $\mys$ and $\myt$ in $\mathcal{M}(n)$ if $\myt$ is formed from $\mys$ by adding a snake of length $i$ (when $i$ is even) or removing a snake of length $i$ (when $i$ is odd). 
From here on, the notation $\mathcal{M}(n)$ will denote this Snakes on a Ming'antu Square puzzle digraph. 
We will construct a diamond-colored distributive lattice $\mathcal{C}(n)$ that is isomorphic to $\mathcal{M}(n)$ as an edge-colored directed graph. 
Let 
\[\mathcal{C}(n) := \left\{\begin{array}{c}\mbox{integer sequences }\\ \mys = (\mys_{1}, \mys_{2},\ldots, \mys_{n})\end{array}\, \rule[-5mm]{0.2mm}{12mm}\, \begin{array}{c} \mys_{1} \geq \mys_{2} \geq \cdots \geq \mys_{n}\\ \mbox{and } 0 \leq \mys_{i} \leq n+1-i \mbox{ for all $i \in [1,n]_{\mathbb{Z}}$}\end{array}\right\},\] 
and partially order these $n$-tuples by component-wise comparison. 
In this partial order, we have a covering relation $\mys \rightarrow \myt$ if and only if there is some $q \in [1,n]_{\mathbb{Z}}$ such that $\mys_{q}+1=\myt_{q}$ but $\mys_{p}=\myt_{p}$ when $p \not= q$; in this case we give the directed edge $\mys \rightarrow \myt$ color $i$ and write $\mys \myarrow{i} \myt$ if $i = n+q-\myt_{q}$. 
Note that each element of $\mathcal{C}(n)$ can be thought of as a partition whose associated partition diagram (aka Ferrers diagram) fits inside an inverted $n \times n$ staircase; see \CatalanLatticeFigure. 
The objects in $\mathcal{C}(n)$ are in evident one-to-one correspondence with the `Dyck lattice paths' that move from the origin at $(0,0)$ to the point $(n+1,n+1)$ in the $xy$-plane taking unit steps in the North and East directions only without landing strictly below the line $y=x$. 
The collection of such lattice paths is well-known to be enumerated by the $(n+1)^{\mbox{\footnotesize st}}$ Catalan number $\frac{1}{n+2}\mbox{\footnotesize $\bigg(\!\!\!\!\begin{array}{c}\mbox{\footnotesize $2(n+1)$}\\ \mbox{\footnotesize $n+1$}\end{array}\!\!\!\!\bigg)$}$, cf.\ \cite{StanleyCatalan}.  
The reader's challenge now is to produce a natural isomorphism $\mathscr{M}: \mathcal{C}(n) \longrightarrow \mathcal{M}(n)$ of edge-colored directed graphs. 

Assuming the existence of such a function, we record without proof the following straightforward analog of \TypeBIsoAndMixedmiddleProps, with definitions of the quantities `$C(\myx,\myy)$' and `$\#_{j}(\myx,\myy)$' modified as follows:  
Let $\myx = (\myx_{1}, \myx_{2},\ldots, \myx_{n})$ and $\myy = (\myy_{1}, \myy_{2},\ldots, \myy_{n})$ be any $n$-tuples in $\mathcal{C}(n)$ with $\myx \leq \myy$. 
Let $C(\myx,\myy)$ denote the set $\left\{j \in [1,n]_{\mathbb{Z}}\  \rule[-1.5mm]{0.2mm}{5mm}\  n+i-\myx_{i} > j \geq n+i-\myy_{i} \mbox{ for some } i \in [1,n]_{\mathbb{Z}}\right\}$; and, for any $j \in [1,n]_{\mathbb{Z}}$, let $\#_{j}(\myx,\myy) := \rule[-1.5mm]{0.2mm}{5mm}\{i \in [1,n]_{\mathbb{Z}}\ |\, n+i-\myx_{i} > j \geq n+i-\myy_{i}\}\rule[-1.5mm]{0.2mm}{5mm}$. 

\begin{figure}[t]
\begin{center}
\hspace*{0.8cm}
\setlength{\unitlength}{1.325cm}
\begin{picture}(6,4.25)
\put(2,0){\VerticalEdgeForBnLattice{Green}{3}{0}{0}}
\put(2,0.6667){\NEEdgeForBnLattice{blue}{2}{0}{0}}
\put(2,0.6667){\NWEdgeForBnLattice{orange}{4}{0}{0}}
\put(1,1.3334){\NWEdgeForBnLattice{Purple}{5}{0}{0}}
\put(1,1.3334){\NEEdgeForBnLattice{blue}{2}{0}{0}}
\put(3,1.3334){\NWEdgeForBnLattice{orange}{4}{0}{0}}
\put(3,1.3334){\NEEdgeForBnLattice{red}{1}{0}{0}}
\put(0,2){\NEEdgeForBnLattice{blue}{2}{0}{0}}
\put(2,2){\NWEdgeForBnLattice{Purple}{5}{0}{0}}
\put(2,2){\VerticalEdgeForBnLattice{Green}{3}{0}{0}}
\put(2,2){\NEEdgeForBnLattice{red}{1}{0}{0}}
\put(4,2){\NWEdgeForBnLattice{orange}{4}{0}{0}}
\put(1,2.6667){\VerticalEdgeForBnLattice{Green}{3}{0}{0}}
\put(1,2.6667){\NEEdgeForBnLattice{red}{1}{-0.2}{-0.2}}
\put(2,2.6667){\NWEdgeForBnLattice{Purple}{5}{-0.2}{0.1}}
\put(2,2.6667){\NEEdgeForBnLattice{red}{1}{0.2}{0.1}}
\put(3,2.6667){\NWEdgeForBnLattice{Purple}{5}{0.2}{-0.2}}
\put(3,2.6667){\VerticalEdgeForBnLattice{Green}{3}{0}{0}}
\put(1,3.3334){\NEEdgeForBnLattice{red}{1}{0}{0}}
\put(2,3.3334){\VerticalEdgeForBnLattice{Green}{3}{0}{0}}
\put(3,3.3334){\NWEdgeForBnLattice{Purple}{5}{0}{0}}
\put(2,0){\CatalanVertexRankZero{-0.4}{0}}
\put(2,0.6667){\CatalanVertexRankOne{0.3}{-0.3}}
\put(1,1.3334){\CatalanVertexRankTwoNumOne{-4.7}{-1}}
\put(3,1.3334){\CatalanVertexRankTwoNumTwo{-0.5}{-1}}
\put(0,2){\CatalanVertexRankThreeNumOne{-4.6}{0}}
\put(2,2){\CatalanVertexRankThreeNumTwo{-2.5}{-2.3}}
\put(4,2){\CatalanVertexRankThreeNumThree{-0.4}{0}}
\put(1,2.6667){\CatalanVertexRankFourNumOne{-4.6}{1}}
\put(2,2.6667){\CatalanVertexRankFourNumTwo{-0.3}{-0.7}}
\put(3,2.6667){\CatalanVertexRankFourNumThree{-0.4}{1}}
\put(1,3.3334){\CatalanVertexRankFiveNumOne{-4}{1}}
\put(2,3.3334){\CatalanVertexRankFiveNumTwo{-4.2}{0.5}}
\put(3,3.3334){\CatalanVertexRankFiveNumThree{-0.4}{1}}
\put(2,4){\CatalanVertexRankSix{-4}{0.7}}
\put(-2.5,-0.75){\small {\bf \CatalanLatticeFigure}\ \  A depiction of $\mathcal{C}(3)$ with its coordinates construed as partitions.}
\end{picture}
\end{center}

\vspace*{0.25cm}
\hspace*{4.25cm}$\overline{\hspace*{7cm}}$
\end{figure}

\noindent
{\bf \CatalanIsoProposition}\ \ 
{\sl Assume there is an isomorphism $\mathscr{M}$ from the diamond-colored distributive lattice $\mathcal{C}(n)$ to the Snakes on a Ming'antu Square puzzle digraph $\mathcal{M}(n)$. 
Then, each diamond-colored distributive lattice $\mathcal{C}(n)$ and $\mathcal{M}(n)$ has $\frac{1}{n+2}\mbox{\footnotesize $\bigg(\!\!\!\!\begin{array}{c}\mbox{\footnotesize $2(n+1)$}\\ \mbox{\footnotesize $n+1$}\end{array}\!\!\!\!\bigg)$}$ elements (the $(n+1)^{\mbox{\footnotesize st}}$ Catalan number) and length $\mbox{\footnotesize $\bigg(\!\!\!\!\begin{array}{c}\mbox{\footnotesize $n+1$}\\ \mbox{\footnotesize $2$}\end{array}\!\!\!\!\bigg)$}$. 
Now let $\mys$ and $\myt$ be objects in $\mathcal{M}(n)$ to be thought of as starting and terminating objects for playing Snakes on a Ming'antu Square. Let $\mys' := \mathscr{M}^{-1}(\mys)$ and $\myt' := \mathscr{M}^{-1}(\myt)$.  Let $\rho_{_\mathcal{C}}$ denote the rank function for $\mathcal{C}(n)$. Then:} {\sl 
\begin{enumerate}
\item[{\sl (1)}] {\sl We have $\rho_{_\mathcal{C}}(\mys')=\sum_{i=1,...,n} \mys'_i$ and $\rho_{_\mathcal{C}}(\myt')=\sum_{i=1,...,n} \myt'_i$.  
Also, $\mys' \vee \myt' = (\mathsf{max}(\mys'_i,\myt'_i))_{i=1,\ldots,n}$ and $\mys' \wedge \myt' = (\mathsf{min}(\mys'_i,\myt'_i))_{i=1,\ldots,n}$. 
That is, the join and meet of  $\mys'$ and $\myt'$ are just the component-wise max and component-wise min, respectively.}  
{\sl Moreover, we have} 
\[\mathsf{dist}(s,t)\ =\ \sum_{i=1}^{n}\left(\rule[-1.5mm]{0mm}{5mm}2\mathsf{max}(\mys'_{i},\myt'_{i}) - \mys'_{i} - \myt'_{i}\rule[-1.5mm]{0mm}{5mm}\right)\ =\ \sum_{i=1}^{n}\left(\rule[-1.5mm]{0mm}{5mm}\mys'_{i} + \myt'_{i} - 2\mathsf{min}(\mys'_{i},\myt'_{i})\rule[-1.5mm]{0mm}{5mm}\right)\ \leq\ \mbox{\footnotesize $\bigg(\!\!\!\begin{array}{c}\mbox{\footnotesize $n+1$}\\ \mbox{\footnotesize $2$}\end{array}\!\!\!\bigg)$}\ <\ \infty.\]
{\sl In particular, the preceding calculation of $\mathsf{dist}(\mys,\myt)$ explicitly solves question (1a) for the Snakes on a Ming'antu Square Puzzle.}
\item[{\sl (2)}] {\sl Keep the notation of part (1), and let $J \subseteq \{1,2,\ldots,n\}$.  
Question (2a) for the Snakes on a Ming'antu Square Puzzle is solved as follows. 
The minimum number of snake moves with colors from the set $J$ that are required in order to obtain $\myt$ from $\mys$ is the left-hand-side quantity in the following inequality:} 
\[\sum_{j \in J}\left(\#_{j}(\mys',\mys'\vee\myt')\rule[-1.5mm]{0mm}{5mm}+\rule[-1.5mm]{0mm}{5mm}\#_{j}(\myt',\mys'\vee\myt')\right)\ \leq\ \mathsf{dist}(\mys,\myt).\] 
{\sl This minimum count of color $J$ snake moves is attained by any shortest sequence of snake moves used to maneuver from $\mys$ to $\myt$. 
Moreover, $\myt$ can be obtained from $\mys$ via snake moves with colors only from the set $J$ if and only if $C(\mys',\mys'\vee\myt') \cup C(\myt',\mys'\vee\myt') \subseteq J$ if and only if} 
\[\mathsf{dist}(\mys,\myt)\ =\ \sum_{j \in J}\left(\#_{j}(\mys',\mys'\vee\myt')\rule[-1.5mm]{0mm}{5mm}+\rule[-1.5mm]{0mm}{5mm}\#_{j}(\myt',\mys'\vee\myt')\right).\]
{\sl All the statements of part (2) remain valid if we appropriately replace each occurrence of ``$\mys'\vee\myt'$'' with ``$\mys'\wedge\myt'$''.}
\end{enumerate}
}

%
\noindent 
\vspace*{-0.5in}
\renewcommand{\refname}{\bf \Large References}
\renewcommand{\baselinestretch}{1.1}

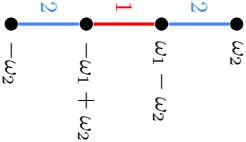
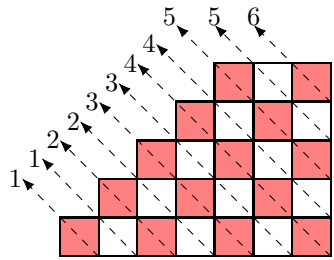
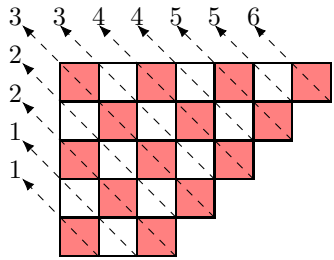
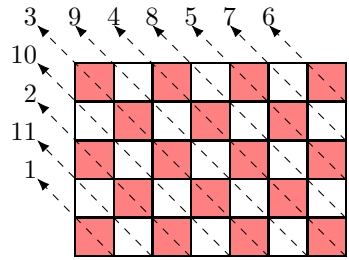
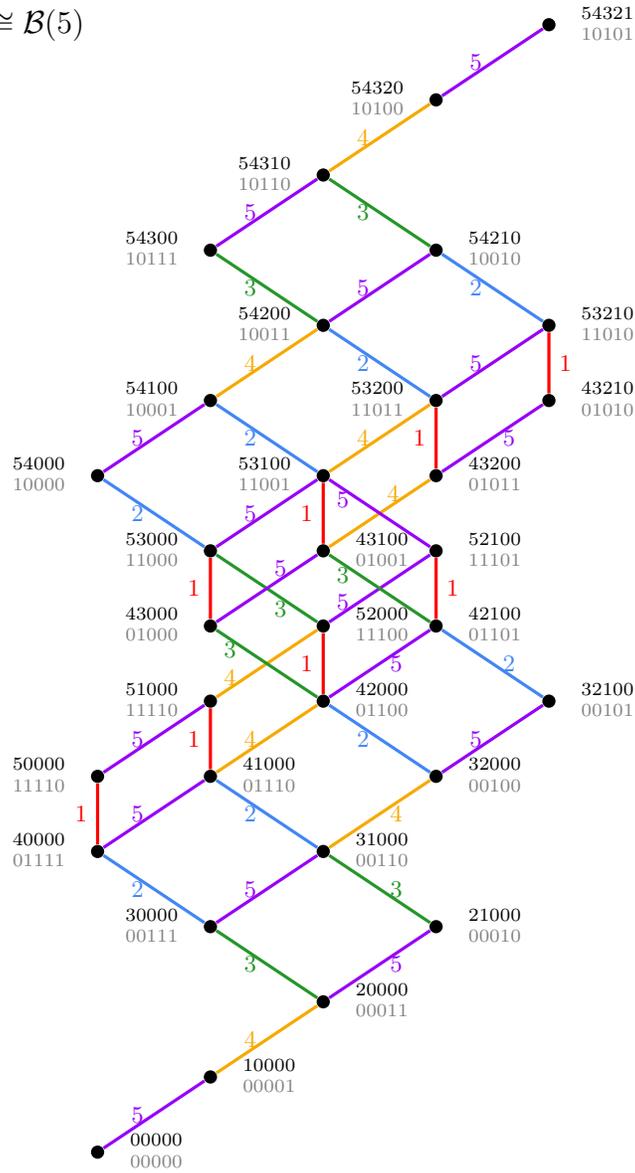
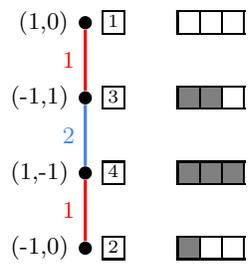
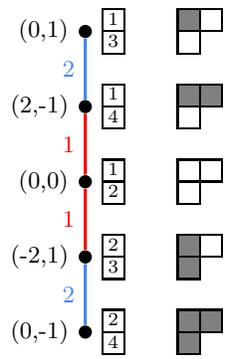
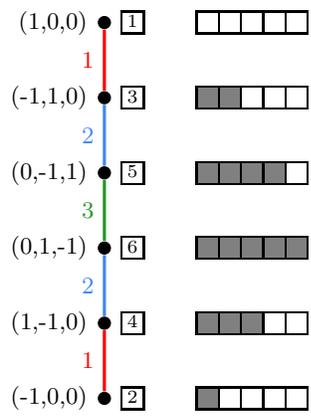
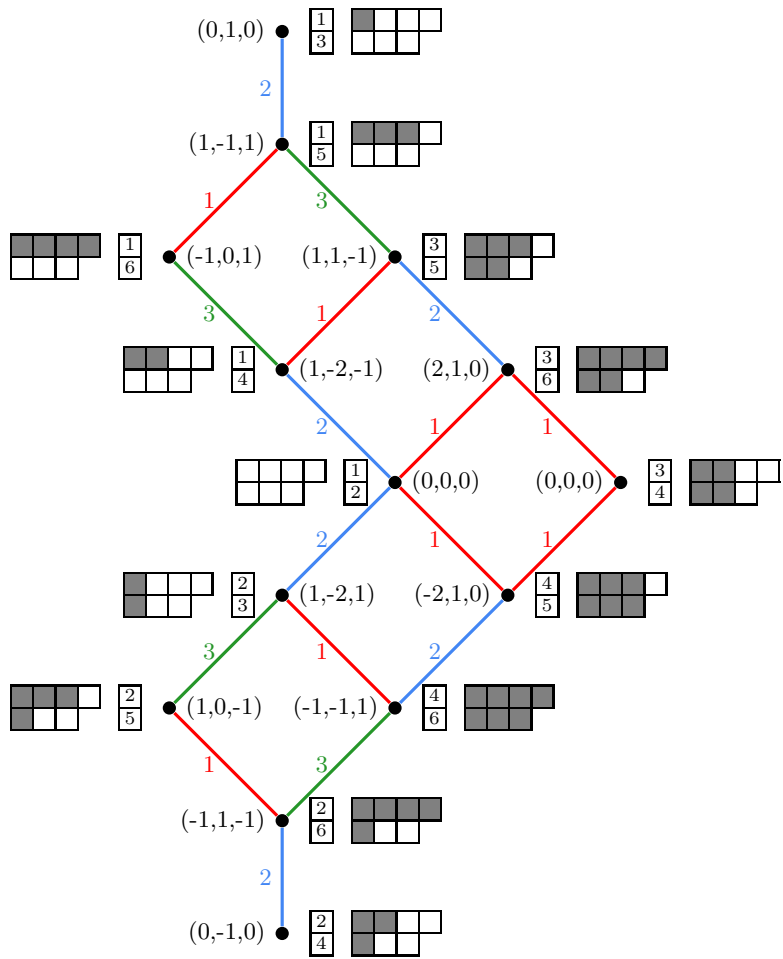
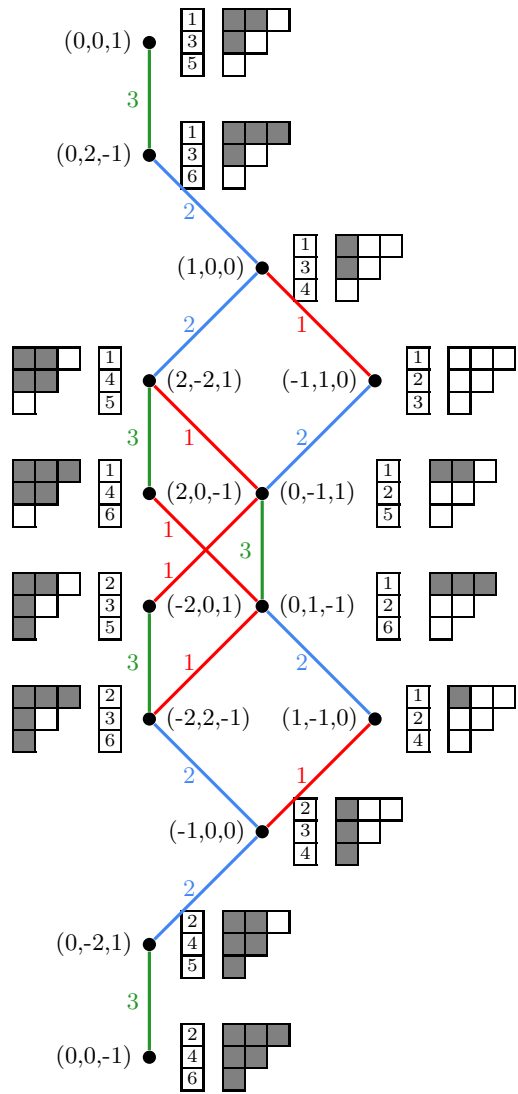
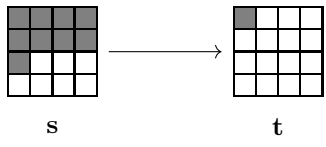
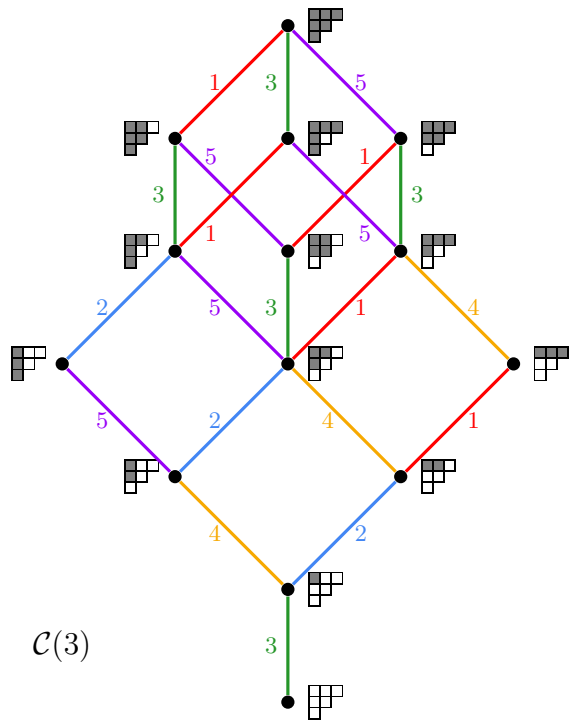
\end{document}